\documentclass[french,a4paper,11pt]{article}

\usepackage[francais]{babel}
\usepackage{eucal}
\usepackage[latin1]{inputenc}
\usepackage[T1]{fontenc}
\usepackage{amsthm}
\usepackage{amsmath,amssymb}
\usepackage{a4wide}
\usepackage{enumerate}
\usepackage{pstricks}
\usepackage{mathrsfs}
\usepackage{eufrak}
\usepackage{makeidx}

\everymath{\displaystyle}

\newtheorem{defn}{Définition}[subsection]
\newtheorem{thm}[defn]{Théorème}
\newtheorem{prop}[defn]{Proposition}
 \newtheorem{cor}[defn]{Corollaire}
 \newtheorem{lem}[defn]{Lemme}
 \newtheorem{rem}[defn]{Remarque}

\newcommand{\C}{\mathbb{C}}

 \newcommand{\Z}{\mathbb{Z}}
 \newcommand{\N}{\mathbb{N}}
 \newcommand{\R}{\mathbb{R}}

 \newcommand{\M}{\mathcal{M}}

\newcommand{\D}{\mathscr{D}}

 \newcommand{\F}{\mathscr{F}}

 \renewcommand{\leq}{\leqslant}
 \renewcommand{\geq}{\geqslant}

\DeclareMathAlphabet{\mathpzc}{OT1}{pzc}{m}{it}

\def\Spin{{\mathrm{Spin}}}
\def\CSpin{{\mathrm{CSpin}}}
\def\CO{{\mathrm{CO}}}
\def\SO{{\mathrm{SO}}}
\def\O{{\mathrm{O}}}
\def\Ric{{\mathrm{Ric}}}
\def\Rd{{\mathrm{R}}}
\def\tr{{\mathrm{tr}}}

\newcommand{\tld}[1]{\widetilde{#1}}
\newcommand{\tra}[1]{\,{\vphantom{#1}}^{t}\!{#1}}

\usepackage[all]{xy}

\begin{document}

\title{Structures conformes asymptotiquement plates}
\author{Guillaume Vassal}
\date{\today}

\maketitle

\begin{abstract}
Dans la première partie de cet article, nous présentons la théorie des structures de Weyl et des spineurs à poids sur une variété conforme. Dans la seconde partie, nous introduisons la notion de structure de Weyl asymptotiquement plate, la masse associée et nous démontrons le théorème de la masse positive dans le cadre de la géométrie spinorielle conforme.
\vskip .6cm


\medskip
\noindent{\it mots clés :} Structure de Weyl, spineur à poids, variété asymptotiquement plate, théorème de la masse positive.
\end{abstract}

\tableofcontents

\section*{Introduction}

Une variété riemannienne $M$ est asymptotiquement plate s'il existe un compact $K$ tel que les composantes connexes de $M\setminus K$ sont difféomorphes à l'espace $\R^{n}$ privé d'une boule, et si dans ces ouverts, la métrique est asymptotique à la métrique plate de $\R^{n}$. Sous certaines conditions, un invariant géométrique calculé à l'infini leur est associée : c'est la \emph{masse} de la variété asymptotiquement plate. Dans un cadre spinoriel, E. Witten \cite{wit} démontre que cette masse est positive si la courbure scalaire de la variété est positive et que ces variétés possèdent une certaine rigidité puisque la masse est nulle si et seulement si la variété est isométrique à l'espace euclidien. Rendue rigoureuse par T. Parker et C. H. Taubes \cite{park}, la méthode repose sur la formule de Lichnerowicz. Une fois intégrée, cette formule fournit un terme de bord qui, sous les bonnes hypothèses, converge vers la masse de la variété riemannienne asymptotiquement plate. Cette formule permet alors de montrer le résultat de positivité de la masse, ainsi que le cas d'égalité.

P. T. Chru{\'s}ciel et M. Herzlich \cite{ch} ont défini une masse pour des variétés asymptotiquement hyperboliques et ont démontré un théorème de masse positive dans ce cadre. Dans \cite{xd}, Xianzhe Dai généralise le théorème de masse positive au cas où la variété est asymptotique au produit $\R^{n}\times X$, où $X$ est une variété compacte simplement connexe de Calabi-Yau ou une variété hyper-Kählerienne. Récemment, V. Minerbe \cite{min}, en s'appuyant sur la méthode de E. Witten, a établi le théorème de masse positive pour des variétés complètes non compactes asymptotiques à une fibration en cercles sur une base euclidienne, dont les fibres sont asymptotiquement de longueur constante.

Dans cet article, nous allons nous intéresser au cas d'une variété $M$ munie d'une classe conforme $c$. Nous allons alors généraliser la notion de platitude asymptotique aux structures de Weyl sur $(M,c)$. Nous dirons que $(M,c,D)$ est une structure de Weyl asymptotiquement plate lorsque $D$ est une connexion de Weyl sur $(M,c)$ et s'il existe une métrique $g$ dans la classe conforme $c$ telles que la variété riemannienne $(M,g)$ est asymptotiquement plate et que la $1$-forme de Lee de $D$ relative à $g$ possède un bon comportement asymptotique à l'infini. Nous définirons alors une notion de masse associée à cette structure généralisant la masse d'une variété asymptotiquement plate. Nous verrons que cette masse est invariante pour certaines sous-classes de métriques de la classe conforme et nous étudierons la façon dont elle évolue par passage d'une classe à l'autre.

Dans le cadre de la théorie des spineurs à poids, P. Gauduchon \cite{gau1} et A. Moroianu \cite{mo} ont établi la formule de Lichnerowicz conforme pour les structures de Weyl sur les fibrés de spineurs à poids. La masse des structures de Weyl asymptotiquement plates apparaît naturellement dans l'intégrale du terme de bord de cette formule. Ce calcul nous permet de montrer sa positivité dans le cas où la courbure scalaire de la connexion de Weyl est positive.

Nous commencerons par des rappels de géométrie conforme sur les structures de Weyl, et les spineurs à poids. Nous démontrerons la formule de Lichnerowicz conforme et donnerons l'expression de son terme de bord. Après de brefs rappels dans le cadre riemannien, nous définirons ensuite avec précision la notion de structure conforme asymptotiquement plate et de la masse associée en donnant ses propriétés d'invariance. Enfin, nous établirons sur le modèle de la démonstration de E. Witten un théorème de masse positive pour les structures de Weyl asymptotiquement plates sur des variétés spinorielles.

Ce texte constitue une partie de ma thèse dirigée par Paul Gauduchon et de Andrei Moroianu. Je les remercie vivement pour leurs idées et pour tout le temps qu'ils m'accordent.

\section{Poids conformes et spineurs}\label{spineur à poids}

Sauf mention contraire, les variétés, les métriques et les sections des différents fibrés considérés seront $C^{\infty}$.

\subsection{Structures de Weyl}

Soit $M$ une variété orientée de dimension $n$, munie d'une structure conforme définie positive $c$.
Rappelons quelques faits de géométrie conforme afin de pouvoir définir la notion de structure de Weyl, introduite par Hermann Weyl dans \cite{weyl}.
Nous noterons $TM$ et $T^{\ast}M$ respectivement les fibrés tangent et cotangent de $M$, $Gl(M)$ ($Gl^{+}(M)$) le fibré
des repères (orientés) de $TM$ et $CO(M)$ le sous-fibré de $Gl(M)$ des repères $c$-orthonormés. Le fibré $CO(M)$ que nous appelons fibré des repères conformes est un fibré principal sur $M$ dont le groupe structural est le groupe conforme $\CO(n)=\O(n)\times\R^{>0}$, où $\R^{>0}$ est l'ensemble des réels positifs non nuls. Sur toute variété différentiable $M$ de dimension $n$, nous pouvons définir une famille $L^{(k)}$, pour $k\in\R$, de fibrés en droites réelles sur $M$ par :
$$L^{(k)}=Gl(M)\times_{|\det|^{k/n}}\R$$
Nous dirons que $L^{(k)}$ est le fibré des scalaires de poids $k$. Notons $L$ les scalaires de poids $1$. Pour $k\in\N$, $L^{(k)}$ est la puissance tensorielle $k$-ième de $L$. Ces fibrés sont naturellement orientables, donc triviaux. De plus, lorsque la variété est munie d'une structure conforme, nous avons une réduction au fibré des repère conformes :
$$L^{(k)}=CO(M)\times_{|\det|^{k/n}}\R$$

Si $\nu$ est une représentation du groupe $\CO(n)$ sur un espace vectoriel $V$, la restriction de $\nu$ au sous-groupe des réels strictement positifs $\R^{>0}$ est de la forme :
$$\nu(a)=a^{k}\mathrm{Id}$$
pour $a\in\R^{>0}$ et où $\mathrm{Id}$ est l'identité de $V$. Nous dirons que le réel $k$ est le poids conforme de la représentation $\nu$. Soit $Gl(n)$ ($Gl^{+}(n)$) le groupe des isomorphismes (orientés) de $\R^{n}$. Les fibrés vectoriels sur $(M,c)$ associés à $Gl(M)$ et obtenus par une représentation de $\mathrm{Gl}(n)$ possèdent un poids conforme naturel qui est le poids de la restriction de cette représentation au sous-groupe conforme $\CO(n)$.  Avec cette convention, les fibrés tangent, cotangent et des $p$-formes sur $M$ sont respectivement de poids naturels $+1$, $-1$ et $-p$. Le réel $k$ du fibré $L^{(k)}$ est le poids de l'action de $\R^{>0}$. Il est facile de voir que le poids naturel d'un produit tensoriel est la somme des poids des facteurs constituant ce produit. Par exemple, le poids conforme naturel de $T^{\ast}M\otimes T^{\ast}M$ est $-2$.

Le complémentaire de la section nulle du fibré $L^{(k)}$ possède deux composantes connexes qui sont données par :
$$L^{(k)}_{+}=CO(M)\times_{|\det|^{k/n}}\R^{>0}\,\ \text{et}\,\ L^{(k)}_{-}=CO(M)\times_{|\det|^{k/n}}\R^{<0}$$
Nous utiliserons aussi le fibré complexifié $L_{\C}^{(k)}$ du fibré des scalaires à poids $L^{(k)}$.

\begin{defn}
Une structure de Weyl sur une variété conforme $(M,c)$ est une connexion linéaire sans torsion sur $TM$ induite par une connexion $\CO(n)$-équivariante sur $\CO(M)$.
\end{defn}

Le théorème fondamental de la géométrie conforme de H. Weyl \cite{gau2} est le suivant (\emph{cf.} aussi \cite{weyl} et \cite{foll}) :  

\begin{thm}\label{thm1}
L'application qui, à toute connexion linéaire $D$ sur $TM$ associe la connexion induite $\nabla^{D}$ sur $L$ détermine, par restriction, un isomorphisme affine de l'espace des structures de Weyl sur $TM$ sur l'espace des connections linéaires sur $L$.
\end{thm}

Soient $D$ une structure de Weyl sur $TM$ et $\nabla^{D}$ sa connexion linéaire associée sur $L$. Nous avons la formule de Koszul généralisée suivante :
\begin{eqnarray}
2c(D_{X}Y,Z)&=&\nabla^{D}_{X}\big(c(Y,Z)\big)+\nabla^{D}_{Y}\big(c(Z,X)\big)-\nabla^{D}_{Z}\big(c(Y,X)\big)\nonumber\\
&&+c(Z,[X,Y])-c(Y,[X,Z])-c(X,[Y,Z])\nonumber
\end{eqnarray}
Cette formule démontre l'existence et l'unicité de $D$ étant donné $\nabla^{D}$. 

Soit $D^{\prime}$ une autre structure de Weyl sur $M$, avec $\nabla^{D^{\prime}}$ sa connexion linéaire sur $L$ associée. La différence $\nabla^{D}-\nabla^{D^{\prime}}$ est une $1$-forme $\theta$ sur $TM$. L'isomorphisme entre les structures de Weyl et les connexions linéaires sur $L$ nous donne la relation suivante :
\begin{equation}\label{eq1}
D_{X}Y=D^{\prime}_{X}Y+\theta(X)Y+(\theta\wedge X)(Y)
\end{equation}
L'endomorphisme antisymétrique $\theta\wedge X$ de $TM$ est défini par :
$$(\theta\wedge X)(Y)=\theta(X)Y-g(X,Y)\theta^{\sharp}$$
où $g$ est une métrique quelconque dans la classe conforme, et où $\theta^{\sharp}$ est le dual riemannien de $\theta$ relativement à $g$. En particulier, nous pouvons associer à une métrique riemannienne $g$ dans la  classe conforme sa $1$-forme de Lee, notée $\theta_{g}$, définie par $\theta_{g}=\nabla^{D}-\nabla^{g}$, où $\nabla^{g}$ est la connexion linéaire sur $L$ induite par la connexion de Levi-Civita de $g$. Notons que la connexion de Levi-Civita d'une métrique $g$ dans $c$, notée $D^{g}$, est une structure de Weyl sur $TM$.

\begin{defn}\label{defn2}
Une structure de Weyl $D$ est fermée (respectivement exacte) si sa connexion linéaire associée $\nabla^{D}$ sur le fibré $L$ est plate (respectivement si $L$ admet une section globale $\nabla^{D}$-parallèle). De manière équivalente, une structure de Weyl $D$ est fermée (respectivement exacte) si $D$ est localement (respectivement globalement) la connexion de Levi-Civita d'une métrique dans la classe conforme.
\end{defn}

Rappelons que la classe conforme $c$ est une section du fibré $S^{2}(T^{\ast}M)\otimes L^{(2)}$ telle que, pour tout champ de vecteurs $X$, $c(X,X)$ est une section de $L_{+}^{(2)}$. La famille $\{e_{i}\}_{i=1\ldots n}$ est une base $c$-orthonormée de $TM$ s'il existe une section positive $l$ de $L$ telle que $c(e_{i},e_{j})=\delta_{ij}l^{2}$. En particulier, la base duale $\{e_{i}^{\ast}\}$ de $\{e_{i}\}$ est donnée par :
$$e_{i}^{\ast}(X)=c(e_{i},X)l^{-2}$$
Nous pouvons étendre les isomorphismes musicaux dans le cadre conforme : nous définissons
d'une part $\flat$ : $TM\rightarrow T^{\ast}M\otimes L^{(2)}$ par :
$$X^{\flat}=c(X,\cdot)$$
pour tout vecteur $X$ de $TM$.
D'autre part,  $\sharp$ : $T^{\ast}M\rightarrow TM\otimes L^{(-2)}$ défini par :
$$\alpha=c(\alpha^{\sharp},\cdot)$$
pour toute $1$-forme $\alpha$ sur $M$.
Par conséquent, nous avons l'identification suivante :
$$T^{\ast}M\otimes T^{\ast}M\cong T^{\ast}M\otimes TM\otimes L^{(-2)}$$
Par contraction des formes et des vecteurs, nous obtenons alors une application de $T^{\ast}M\otimes T^{\ast}M$ dans $L^{(-2)}$ : c'est la trace conforme des formes bilinéaires symétriques.

Pour tous $X$, $Y$ et $Z$ vecteurs de $TM$, nous définissons le tenseur de courbure de la structure de Weyl $D$, que nous notons $\Rd^{D}$, par :
$$\Rd^{D}_{X,Y}Z=[D_{X},D_{Y}]Z-D_{[X,Y]}Z$$
Nous pouvons voir $\Rd^{D}$ comme une section du fibré $T^{\ast}M^{\otimes^{3}}\otimes TM$. Nous définissons alors $\Ric^{D}$ le tenseur de Ricci de la structure de Weyl $D$ par :
$$\Ric^{D}(X,Y)=\mathrm{trace}(Z\mapsto\Rd^{D}_{Z,X}Y)$$
pour $X$ et $Y$ vecteurs de $TM$. Le tenseur de Ricci est une section du fibré $T^{\ast}M\otimes T^{\ast}M$. En appliquant la trace conforme sur le tenseur de Ricci, nous obtenons une section de $L^{(-2)}$. Cette section, que nous notons $Scal^{D}$, est la courbure scalaire de la structure de Weyl.

Nous avons une correspondance entre les métriques riemanniennes dans $c$ et les sections du fibré $L_{+}^{2}$ : une métrique riemannienne $g$ et la section $l^{2}$ du fibré $L_{+}^{(2)}$ correspondante sont reliées par :
$$c=g\otimes l^{2}$$
Le choix d'une métrique dans la classe conforme détermine, à un élément du groupe orthogonal $\O(n)$ près, un repère orthonormé au voisinage de chaque point de $M$. Cependant, une section $l^{k}$ du fibré $L^{(k)}$ est une classe d'équivalence sous l'action de $Gl(n)$ de la forme $[s,v]$, où $s$ est un repère local de $TM$ et $v$ une fonction sur $M$. Par conséquent, lorsqu'une métrique est fixée dans la classe conforme, les sections des fibrés $L^{(k)}$ s'identifient à des fonctions sur la variété.
Soit $g$ une métrique dans classe conforme ; la structure de Weyl $D$ et la connexion de Levi-Civita $D^{g}$ sont liées via la $1$-forme de Lee $\theta_{g}$ par la formule suivante :
\begin{equation}\label{eq2}
D_{X}Y=D^{g}_{X}Y+\theta_{g}(X)Y+(\theta_{g}\wedge X)(Y)
\end{equation}
Ainsi, la courbure scalaire de $D$ s'identifie, via la trivialisation de $L$ par $g$, à une fonction sur $M$ par la formule suivante \cite{gau2} :
\begin{equation}\label{eq3}
Scal^{D}=Scal^{g}-2(n-1)\mathrm{tr}_{g}(D^{g}\theta_{g})-(n-1)(n-2)|\theta_{g}|_{g}^{2}
\end{equation}
où $Scal^{g}$ est la courbure scalaire de la connexion de Levi-Civita $D^{g}$ et $\mathrm{tr}_{g}$ la trace relative à la métrique $g$. Pour plus d'information concernant les structures de Weyl, le lecteur intéressé pourra consulter \cite{gau2}.

\subsection{Spineurs conformes}

Soit $M$ une variété spinorielle de dimension $n$ \cite{law}. L'algèbre de Clifford réelle $Cl_{n}$ associée à l'espace $(\R^{n},\|\, \|)$ euclidien est l'unique algèbre réelle, à isomorphisme près, vérifiant la propriété universelle suivante : pour toute $\R$-algèbre associative unitaire $A$, une application linéaire $v$ de $\R^{n}$ dans $A$ telle que $v(x)^{2}=-\|x\|^{2}1_{A}$, pour tout $x$ de $\R^{n}$, s'étend de façon unique en un morphisme d'algèbres de $Cl_{n}$ dans $A$. Soient $\Spin(n)\subset Cl_{n}$ le groupe spinoriel et $\lambda$: $\Spin(n)\rightarrow\SO(n)$ le revêtement à deux feuillets du groupe spécial orthogonal $\SO(n)$. Soit $\mu$ : $\Spin(n)\rightarrow\mathrm{Aut}(\Delta_{n})$ la représentation spinorielle de $\Spin(n)$ sur l'espace des spineurs $\Delta_{n}$. La représentation spinorielle est la restriction de la représentation de l'algèbre de Clifford $Cl_{n}$ sur l'espace $\Delta_{n}$. Cette représentation induit une action de $Cl_{n}$ sur $\Delta_{n}$ (cette action est aussi appelé multiplication de Clifford). En particulier, nous avons une inclusion canonique de $\R^{n}$ dans $Cl_{n}$, donc une action de $\R^{n}$ sur $\Delta_{n}$. Nous notons $x\cdot\xi$ cet action de Clifford, où $x\in Cl_{n}$ et $\xi\in\Delta_{n}$. Rappelons qu'il existe un isomorphisme d'espace vectoriel entre $Cl_{n}$ et l'algèbre extérieure $\Lambda^{\ast}\R^{n}$ de $\R^{n}$ donné par :
\begin{eqnarray}
\Lambda^{\ast}\R^{n}&\rightarrow& Cl_{n}\nonumber\\
e_{i_{1}}\wedge\ldots\wedge e_{i_{k}}&\mapsto& e_{i_{1}}\cdots e_{i_{k}}\nonumber
\end{eqnarray}
où $(e_{1},\ldots,e_{n})$ est une base de $\R^{n}$. Notons que $v\wedge$ et $v\lrcorner$ sont respectivement le produit extérieur et intérieur par $v$ sur $M$. Nous avons pour l'action de Clifford l'identification suivante :
$$x\cdot(\omega\cdot\xi)=(x\wedge\omega)\cdot\xi-x\lrcorner\omega\cdot\xi$$
où $x\in\R^{n}$, $\omega\in\Lambda^{\ast}\R^{n}$ et $\xi\in\Delta_{n}$. La multiplication de Clifford est un morphisme de $\Spin(n)$-représentations, $i.e.$ :
$$\mu(\gamma)(x\cdot\xi)=(\lambda(\gamma)x)\cdot(\mu(\gamma)\xi)$$
pour $\gamma\in\Spin(n)$, $x\in\R^{n}$ et $\xi\in\Delta_{n}$. Soit $g$ une métrique riemannienne sur $M$. Une structure $spin$ sur $(M,g)$ étant choisie, notons $Spin_{g}(M)$ le $\Spin(n)$-fibré principal au-dessus du fibré des repères $g$-orthonormés directs $SO_{g}(M)$. Le fibré des spineurs relatif à la métrique $g$ est défini par : 
$$\Sigma^{g}=Spin_{g}(M)\times_{\mu}\Delta_{n}$$
Pour plus de détails concernant la géométrie spinorielle dans le cadre riemannien, le lecteur pourra consulter \cite{law}.

Nous définissons le groupe spinoriel conforme $\CSpin(n)$ comme le produit $\Spin(n)\times\R^{>0}$. Pour tout $\tld{\gamma}\in\CSpin(n)$, nous écrirons :
$$\tld{\gamma}=a\gamma$$
où $a\in\R_{+}$ et $\gamma\in\Spin(n)$ sont uniquement déterminés. Soit $\CO^{+}(n)$; nous avons le morphisme de groupes $\tld{\lambda}$ : $\CSpin(n)\rightarrow\CO^{+}(n)$ défini comme le produit de $\lambda$ et de l'identité de $\R^{>0}$ et où $\CO^{+}(n)=\SO(n)\times\R^{>0}$. Soit $k\in\R$. La représentation spinorielle conforme de poids $k$, notée $\mu^{(k)}$, est la représentation linéaire de $\CSpin(n)$ sur l'espace des spineurs $\Delta_{n}$ définie par :
$$\mu^{(k)}(\tld{\gamma})=a^{k}\mu(\gamma)$$
pour tout $\tld{\gamma}$ dans $\CSpin(n)$. De la même façon que la représentation spinorielle, nous avons une relation de compatibilité entre la multiplication de Clifford et la représentation spinorielle conforme de poids $k$.
Pour $\tld{\gamma}\in\CSpin(n)$, $x\in\R^{n}$ et $\xi\in\Delta_{n}$, nous avons :
$$\mu^{(k)}(\tld{\gamma})(x\cdot\xi)=\left(\tld{\lambda}(\tld{\gamma})x\right)\cdot\left(\mu^{(k)}(\tld{\gamma})\xi\right)$$
Le fibré vectoriel des spineurs de poids conforme $k$ est défini par:
$$\Sigma^{(k)}=CSpin(M)\times_{\mu^{(k)}}\Delta_{n}$$
Nous avons l'identification suivante :
$$\Sigma^{(k)}\cong\Sigma^{(0)}\otimes L^{(k)}$$
L'action de Clifford de $TM$ sur l'espace des spineurs à poids est l'application
\begin{eqnarray}
&TM\otimes\Sigma^{(k)}&\rightarrow\Sigma^{(k+1)}\nonumber\\
&X\otimes\psi&\mapsto  X\cdot\psi\nonumber
\end{eqnarray}
définie par :
$$X\cdot\psi=[\tld{s},x\cdot\xi]$$
où $\psi$ et $X$ sont respectivement représentés par $[\tld{s},\xi]$ et $[s,x]$, avec $\tld{s}$ un repère de $\CSpin(M)$, $s=\tld{\lambda}(\tld{s})$ et $x\in\R^{n}$.

\begin{prop}\label{prop1}
Pour tout $g\in c$, et tout $k\in\R$, il existe un isomorphisme canonique de $\Sigma^{g}$ dans $\Sigma^{(k)}$.
\end{prop}

\begin{proof}
Une métrique $g$ dans la classe conforme définit une réduction du fibré $CSpin(M)$ à $Spin_{g}(M)$. Dans ce cas, nous avons :
$$CSpin(M)\times_{\mu^{(k)}}\Delta_{n}=Spin_{g}(M)\times_{\mu}\Delta_{n}$$
puisque $\mu$ est le restriction de la représentation $\mu^{(k)}$ au groupe $\Spin(n)$.
\end{proof}

Soient $g\in c$ et $\widetilde{g}=f^{-2}g$, avec $f$ une fonction réelle non nulle sur $M$. La proposition précédente nous donne une famille d'isomorphismes
\begin{eqnarray}
\Phi^{(k)}\ :&\Sigma^{g}M&\rightarrow\Sigma^{\widetilde{g}}M\nonumber\\
&[\tld{s},v]&\mapsto [\tld{s}f,f^{-k}v]\nonumber
\end{eqnarray}
pour tout $k\in\R$. Soit $\langle\,,\,\rangle$ le produit scalaire hermitien sur $\Delta_{n}$ compatible avec l'action de Clifford. Notons $(\,,\,)_{g}$ le produit scalaire hermitien $\Spin(n)$-invariant sur $\Sigma^{g}$ induit par $\langle\,,\,\rangle$. Les isomorphismes $\Phi^{(k)}$ ne sont pas des isométries :
$$|\Phi^{(k)}\psi|_{\widetilde{g}}=f^{-2k}|\psi|_{g}$$
Nous définissons une application bilinéaire $h$ par :
\begin{eqnarray}
h\ :&\Sigma^{(k_{1})}\otimes\Sigma^{(k_{2})}&\rightarrow L_{\C}^{(k_{1}+k_{2})}\nonumber\\
&\psi\otimes\varphi&\mapsto[s,\langle u,v\rangle]\nonumber
\end{eqnarray}
où $\psi$ et $\varphi$ sont représentés respectivement par $[\tld{s},\xi]$ et $[\tld{s},\zeta]$ et où $s$ est le projeté de $\tld{s}$ sur $CO(M)$.

\begin{prop}\label{prop2}
L'application bilinéaire $h$ est bien définie et pour tous $X\in TM$, $\psi\in\Sigma^{(k_{1})}$ et $\varphi\in\Sigma^{(k_{2})}$, nous avons :
$$h(X\cdot\psi,\varphi)=-h(\psi,X\cdot\varphi)\in L_{\C}^{(k_{1}+k_{2}-1)}$$
\end{prop}

\begin{proof}
En changeant le repère $\tld{s}$ par $\tld{s}\cdot\tld{\gamma}^{-1}$, où $\tld{\gamma}\in\CSpin(n)$, nous obtenons
\begin{eqnarray}
\langle\mu^{(k_{1})}(\tld{\gamma})\xi,\mu^{(k_{2})}(\tld{\gamma})\zeta\rangle
&=&\langle a^{k_{1}}\mu(\gamma)\xi,a^{k_{2}}\mu(\gamma)\zeta\rangle\nonumber\\
&=&a^{k_{1}+k_{2}}\langle \mu(\gamma)\xi,\mu(\gamma)\zeta\rangle\nonumber\\
&=&a^{k_{1}+k_{2}}\langle\xi,\zeta\rangle\nonumber\\
&=&(\det\tld{\gamma})^{(k_{1}+k_{2})/n}\langle\xi,\zeta\rangle\nonumber
\end{eqnarray}
Ce calcul montre que $[s,\langle\xi_{1},\xi_{2}\rangle]$ est une section de $L_{\C}^{(k_{1}+k_{2})}$ et donc que $h$ est bien définie. Pour $\xi$, $\zeta\in\Delta_{n}$ et $x\in\R^{n}$, le produit scalaire hermitien sur $\Delta_{n}$ vérifie :
$$\langle x\cdot\xi,\zeta\rangle=-\langle\xi,x\cdot\zeta\rangle$$
Par définition de la multiplication de Clifford des spineurs à poids, il est clair que $h$ vérifie la même relation.
\end{proof}

Les sections du fibré $L_{\C}^{(-n)}$ sont des densités d'intégration sur $M$. Pour $k_{1}+k_{2}=-n$, nous définissons une application sesquilinéaire $H$ : $C_{0}^{\infty}(\Sigma^{(k_{1})})\times C_{0}^{\infty}(\Sigma^{(k_{2})})\rightarrow\C$ par :
$$H(\psi,\varphi)=\int_{M}h(\psi,\varphi)$$
où $C^{\infty}_{0}$ désigne l'espace des sections lisses à support compact. Soient $k_{1}$ et $k_{2}$ des réels tels que $k_{1}+k_{2}=-n$ et $g$ une métrique dans la classe conforme $c$. Nous avons vu que les fibrés $\Sigma^{(k_{1})}$ et $\Sigma^{(k_{2})}$ s'identifient canoniquement au fibré $\Sigma^{g}$. La multiplication de Clifford définie sur les spineurs à poids s'identifie alors à la multiplication de Clifford (usuelle) sur $\Sigma^{g}$. D'autre part, l'image de l'application $h$ par cette identification est le produit scalaire hermitien $(\,,\,)_{g}$ défini sur $\Sigma^{g}$. Ainsi, nous pouvons remarquer que, pour une métrique $g$ quelconque dans $c$, nous avons :
$$H(\psi,\varphi)=\int_{M}(\psi,\varphi)_{g}v_{g}$$
où $v_{g}$ est la forme volume associée à la métrique $g$.

\begin{prop}\label{prop3}
Soient $D$ et $\widetilde{D}$ deux structures de Weyl sur $(M,c)$ reliées par (\ref{eq1}). Ces dernières induisent deux connexions linéaires $D^{(k)}$ et $\widetilde{D}^{(k)}$ sur $\Sigma^{(k)}$ reliées par :
\begin{equation}\label{eq3bis}
D^{(k)}_{X}\psi=\widetilde{D}^{(k)}_{X}\psi-\frac{1}{2}X\cdot\theta\cdot\psi+(k-\frac{1}{2})\theta(X)\psi\nonumber
\end{equation}
pour tout $\psi$ dans $\Sigma^{(k)}$.
\end{prop}

\begin{proof}
Nous avons vu que les structures de Weyl $D$ et $\widetilde{D}$ sont reliées par $D=\tld{D}+\Theta$, où $\Theta(X)=\theta\wedge X+\theta(X)\mathrm{Id}$. Les connexions induites par $D$ et $\widetilde{D}$ sur le fibré $\Sigma^{(k)}$ vérifient donc la relation suivante :  $$D^{(k)}=\tld{D}^{(k)}+d\mu^{(k)}(\Theta)$$
où $d\mu^{(k)}$ est la différentielle à l'origine de la représentation linéaire de poids $k$ de $\CSpin(n)$. Nous obtenons :
$$D^{(k)}=\tld{D}^{(k)}+\sum_{i=1}^{n}d\mu(\theta\wedge e_{i})\otimes e_{i}^{\ast}+k\mathrm{Id}\otimes\theta$$
Cela donne, pour $\psi\in\Sigma^{(k)}$ et $X\in TM$, la formule suivante :
$$D^{(k)}_{X}\psi=\tld{D}^{(k)}_{X}\psi+\frac{1}{2}(\theta\wedge X)\cdot\psi+k\theta(X)\psi$$
De plus, d'après la définition du produit de Clifford des $2$-formes, nous avons :
$$X\cdot\theta=-\theta\wedge X-\theta(X)$$
Nous obtenons alors la formule souhaitée.
\end{proof}

Soient $D$ une structure de Weyl sur $(M,c)$, $k$ un réel et $D^{(k)}$ la connexion induite par $D$ sur $\Sigma^{(k)}$. Soit $g\in c$. La connexion de Levi-Civita $D^{g}$ de $g$ est une structure de Weyl sur $TM$. La connexion $D^{g}$ détermine donc une connexion linéaire sur le fibré des spineurs à poids $\Sigma^{(k)}$. Nous continuons de noter $D^{g}$ cette connexion. La proposition précédente nous donne la relation suivante :
\begin{equation}\label{eq3ter}
D^{(k)}_{X}\psi=D^{g}_{X}\psi-\frac{1}{2}X\cdot\theta_{g}\cdot\psi+(k-\frac{1}{2})\theta_{g}(X)\psi
\end{equation}
pour tout $\psi$ dans $\Sigma^{(k)}$ et tout $X$ dans $TM$.

Soit $D^{\prime}$ une autre structure de Weyl. Nous avons $\nabla^{D}$ et $\nabla^{D^{\prime}}$ les connexions linéaires sur $L$ associées respectivement aux structures de Weyl $D$ et $D^{\prime}$. Ces connexions induisent des connexions linéaires sur les fibrés $L^{(k)}$ pour $k\in\R$, que l'on note toujours $\nabla^{D}$ et $\nabla^{D^{\prime}}$. Nous avons vu qu'il existe une $1$-forme $\theta$ sur $M$ telle que $\nabla^{D}=\nabla^{D^{\prime}}+\theta$. Par conséquent, pour toute section $l^{k}$ de $L^{(k)}$ et tout vecteur $X$ de $TM$, nous avons : 
$$\nabla^{D}_{X}l^{k}=\nabla^{D^{\prime}}_{X}l^{k}+k\theta(X)l^{k}$$
En particulier, si nous fixons une métrique $g$ de la classe conforme, les sections de $L^{(k)}$ s'identifient à des fonctions sur $M$. Ainsi, lorsque $D^{\prime}$ est la connexion de Levi-Civita de $g$, nous obtenons :
$$\nabla^{D}_{X}l^{k}=dl^{k}(X)+k\theta_{g}(X)l^{k}$$

Soit $\D^{g}$ : $\Sigma^{g}\rightarrow\Sigma^{g}$ l'opérateur de Dirac agissant sur le fibré des spineurs relatif à la métrique $g$. L'opérateur $\D^{g}$ est défini comme la composition de la multiplication de Clifford et de la connexion induite par $D^{g}$ sur $\Sigma^{g}$. Cet opérateur de Dirac est un opérateur différentiel linéaire d'ordre $1$ et elliptique. Pour plus de détails, nous renvoyons de nouveau le lecteur à \cite{law}. La multiplication de Clifford sur les spineurs à poids définit une contraction, notée $m^{(k)}$, sur les spineurs de poids $k$ :
\begin{eqnarray}
m^{(k)}\ :&T^{\ast}M\otimes\Sigma^{(k)}&\rightarrow\Sigma^{(k-1)}\nonumber\\
&\omega\otimes\psi&\mapsto \omega\cdot\psi\nonumber
\end{eqnarray}
Nous considérons la connexion $D^{(k)}$ comme un opérateur de $\Sigma^{(k)}$ dans $T^{\ast}M\otimes\Sigma^{(k)}$. Nous définissons alors l'opérateur de Dirac de poids $k$ par :
$$\D^{(k)}\ :\ \Sigma^{(k)}\xrightarrow[]{m^{(k)}\circ D^{(k)}}\Sigma^{(k-1)}$$
De cette façon, la connexion de Levi-Civita $D^{g}$ agissant sur les sections de $\Sigma^{(k)}$ détermine un opérateur de Dirac de poids $k$. Lorsque les spineurs de poids $k$ sont identifiés à des spineurs de $\Sigma^{g}$, via le choix de la métrique $g$, cet opérateur de Dirac à poids s'identifie tautologiquement à l'opérateur de Dirac usuel $\D^{g}$ sur $\Sigma^{g}$. Les calculs qui suivent sont effectués dans les différents espaces à poids mentionnés. En revanche, si une métrique est fixée, les calculs se ramènent, via les identifications induites par la métrique, au cas riemannien, en identifiant les différents objets à leurs correspondants riemanniens respectifs.

\begin{prop}\label{prop4}
Soient $g$ une métrique de la classe conforme $c$ et $D$ une structure de Weyl sur $(M,c)$ dont la forme de Lee par rapport à $g$ est $\theta_{g}$. Les opérateurs de Dirac $\D^{g}$ et $\D^{(k)}$ agissant sur les sections de $\Sigma^{(k)}$ sont reliés par la formule suivante :
\begin{equation}\label{eq4}
\D^{(k)}\psi=\D^{g}\psi+(k+\frac{1}{2}(n-1))\theta_{g}\cdot\psi\nonumber
\end{equation}
pour toute section $\psi$ de $\Sigma^{(k)}$.
\end{prop}

\begin{proof}
La proposition \ref{prop3} nous donne :
$$D^{(k)}_{X}\psi=D^{g}_{X}\psi+(k-\frac{1}{2})\theta_{g}(X)\psi-\frac{1}{2}X\cdot\theta_{g}\cdot\psi$$
Par conséquent, dans une base $\{e_{i}\}$ $g$-orthonormée, nous avons :
\begin{eqnarray}
\D^{(k)}\psi&=&\D^{g}\psi+(k-\frac{1}{2})\sum_{i=1}^{n}\theta_{g}(e_{i})e_{i}\cdot\psi-\frac{1}{2}\sum_{i=1}^{n}e_{i}\cdot e_{i}\cdot\theta_{g}\cdot\psi\nonumber\\
&=&\D^{g}\psi+(k-\frac{1}{2})\theta_{g}\cdot\psi+\frac{n}{2}\theta_{g}\cdot\psi\nonumber
\end{eqnarray}
\end{proof}

\subsection{Formule de Lichnerowicz conforme I}

Dans le cas des spineurs associés à une métrique riemannienne, la formule de Lichnerowicz relie le carré de l'opérateur de Dirac, l'opérateur de Laplace-Beltrami et la courbure scalaire de la métrique. Fixons une métrique $g$ dans la classe conforme $c$ et une structure de Weyl $D$ sur $M$. Nous avons la proposition suivante \cite{law} : 

\begin{prop}\label{prop5}
Soit $(M,g)$ une variété riemannienne de dimension $n$, orientée, et spinorielle. Pour $\psi$ dans $\Sigma^{g}$, nous avons la formule de Lichnerowicz suivante :
$$(\D^{g})^{2}\psi=\mathrm{tr}(D^{g,2})(\psi)+\frac{1}{4}Scal^{g}\psi$$
où $D^{g,2}$ est la dérivée covariante seconde et $Scal^{g}$ la courbure scalaire de la connexion de Levi-Civita.
\end{prop}

Dans le cadre des spineurs à poids, P. Gauduchon \cite{gau1} et A. Moroianu \cite{mo} ont établi une formule de Lichnerowicz conforme. Soit $k$ un réel. Etant donnée une structure de Weyl, l'opérateur de Dirac de poids $k$ agit sur les sections de spineurs de poids $k$ mais celui-ci fournit des sections de spineurs de poids $k-1$. Par conséquent, nous définissons le carré de l'opérateur de Dirac de poids $k$ comme la composition $\D^{(k-1)}\circ\D^{(k)}$, qui est un opérateur défini sur les sections de $\Sigma^{(k)}$ et à valeurs dans l'espace des sections de $\Sigma^{(k-2)}$.

\begin{prop}\label{prop6}
Soit $\psi$ un spineur de poids $k$. Le carré de l'opérateur de Dirac à poids associé à une structure de Weyl dont la forme de Lee est $\theta_{g}$ est donné par :
\begin{eqnarray}
\D^{2}\psi
&=&\D^{(k-1)}\D^{(k)}\psi\nonumber\\
&=&(\D^{g})^{2}\psi+(k+\frac{1}{2}(n-1))\D^{g}(\theta_{g})\cdot\psi-\theta_{g}\cdot\D^{g}\psi\nonumber\\
&&-(2k+n-1)D^{g}_{\theta_{g}^{\sharp}}\psi-\Big(k^{2}-k(2-n)+\frac{n^{2}-4n+3}{4}\Big)|\theta_{g}|_{g}^{2}\psi\nonumber
\end{eqnarray}
\end{prop}

\begin{proof}
Soient $x$ dans $M$ et $\{e_{i}\}$ une base $g$-orthonormée de $T_{x}M$. Dans cette base, nous avons :
\begin{eqnarray}
&&\D^{(k-1)}\D^{(k)}\psi=\D^{(k-1)}(\D^{g}\psi+(k+\frac{1}{2}(n-1))\theta_{g}\cdot\psi)\nonumber\\
&=&\D^{g}(\D^{g}\psi+(k+\frac{1}{2}(n-1))\theta_{g}\cdot\psi)+(k-1+\frac{1}{2}(n-1))\theta_{g}\cdot(\D^{g}\psi+(k+\frac{1}{2}(n-1))\theta_{g}\cdot\psi)\nonumber\\
&=&(\D^{g})^{2}\psi+(k+\frac{1}{2}(n-1))\D^{g}(\theta_{g})\cdot\psi-(k+\frac{1}{2}(n-1))\theta_{g}\cdot\D^{g}\psi-(2k+(n-1))D_{\theta_{g}^{\sharp}}^{g}\psi\nonumber\\
&&+(k-1+\frac{1}{2}(n-1))\theta_{g}\cdot\D^{g}\psi-(k-1+\frac{1}{2}(n-1))(k+\frac{1}{2}(n-1))|\theta_{g}|_{g}^{2}\psi\nonumber\\&=&(\D^{g})^{2}\psi+(k+\frac{1}{2}(n-1))\D^{g}(\theta_{g})\cdot\psi-\theta_{g}\cdot\D^{g}\psi-(2k+n-1)D^{g}_{\theta_{g}^{\sharp}}\psi\nonumber\\
&&-(k+\frac{1}{2}(n-1))(k-1+\frac{1}{2}(n-1))|\theta_{g}|_{g}^{2}\psi\nonumber
\end{eqnarray}
\end{proof}

Comme pour toute connexion sans torsion, nous définissons la connexion de Weyl seconde de poids $k$ par :
$$D_{X,Y}^{(k),2}\psi=D_{X}^{(k)}(D_{Y}^{(k)}\psi)-D_{D_{X}Y}^{(k)}\psi$$
où $X$ et $Y$ sont des champs de vecteurs de $TM$, et $\psi\in\Sigma^{(k)}$. La connexion de Weyl seconde de poids $k$ peut être vue comme un opérateur agissant sur les sections de $\Sigma^{(k)}$ à valeurs dans l'espace des sections de $T^{\ast}M\otimes T^{\ast}M\otimes\Sigma^{(k)}$. En appliquant la trace conforme et en rappelant que $L^{(-2)}\otimes\Sigma^{(k)}\cong\Sigma^{(k-2)}$, nous obtenons  l'application  $\tr(D^{(k),2})$ : $\Sigma^{(k)}\rightarrow\Sigma^{(k-2)}$.

\begin{prop}\label{prop7}
Soient $\{e_{i}\}$ un repère local $c$-orthonormé tel que $c(e_{i},e_{j})=\delta_{ij}l^{2}$, où $l\in L$, et $\psi$ une section de $\Sigma^{(k)}$. Une écriture locale de $\mathrm{tr}(D^{(k),2})$ est donnée par :
$$\mathrm{tr}(D^{(k),2})(\psi)=-\sum_{i=1}^{n}\left(D_{e_{i}}^{(k)}(D_{e_{i}}^{(k)}\psi)-D^{(k)}_{D_{e_{i}}e_{i}}\psi\right)l^{-2}$$
\end{prop}

\begin{proof}
Montrons que l'expression ne dépend pas du choix du repère $c$-orthonormé. Soit $\{\tld{e}_{i}\}$ une base $c$-orthonormée. Il existe un champ de matrices conformes $A=(a_{ij})$ tel que :
$$\tld{e}_{i}=\sum_{k=1}^{n}a_{ik}e_{k}$$
Nous avons $\tra{A}A=f^{2}$, où $f$ est une fonction sur $M$. De manière équivalente, nous avons la relation suivante :
$$\sum_{k=1}^{n}a_{ik}a_{jk}=\delta_{ij}f^{2}$$
Nous obtenons alors :
\begin{eqnarray}
c(\tld{e}_{i},\tld{e}_{j})&=&\sum_{k=1}^{n}\sum_{l=1}^{n}a_{ik}a_{jl}c(e_{k},c_{l})\nonumber\\
&=&\sum_{k=1}^{n}a_{ik}a_{jk}l^{2}\nonumber\\
&=&\delta_{ij}(fl)^{2}\nonumber
\end{eqnarray}
Par conséquent, nous avons $c(\tld{e}_{i},\tld{e}_{j})=\tilde{l}^{2}$, où $\tld{l}=fl$. Nous calculons :
\begin{eqnarray}
\sum_{i=1}^{n}\left(D^{(k)}_{\tld{e_{i}}}(D^{(k)}_{\tld{e_{i}}}\psi)-D^{(k)}_{D_{\tld{e_{i}}}\tld{e_{i}}}\psi\right)\tld{l}^{-2}
&=&\sum_{i,p,q=1}^{n}\left(a_{ip}a_{iq}D^{(k)}_{e_{p}}(D^{(k)}_{e_{q}}\psi)+D^{(k)}_{\tld{e}_{i}}(a_{iq})D^{(k)}_{e_{q}}\psi\right)\tld{l}^{-2}\nonumber\\
&&-\sum_{i,p,q=1}^{n}\left(a_{ip}a_{iq}D^{(k)}_{D_{e_{p}}e_{q}}\psi+D^{(k)}_{\tld{e}_{i}}(a_{iq})D^{(k)}_{e_{q}}\psi\right)\tld{l}^{-2}\nonumber\\
&=&\sum_{p,q=1}^{n}\Big(\sum_{i=1}^{n}a_{ip}a_{iq}\Big)\left(D^{(k)}_{e_{p}}(D^{(k)}_{e_{q}}\psi)-D^{(k)}_{D_{e_{p}}e_{q}}\psi\right)\tld{l}^{-2}\nonumber\\
&=&\sum_{p,q=1}^{n}\delta_{pq}\left(D^{(k)}_{e_{p}}(D^{(k)}_{e_{q}}\psi)-D^{(k)}_{D_{e_{p}}e_{q}}\psi\right)f^{2}\tld{l}^{-2}\nonumber\\
&=&\sum_{p=1}^{n}\left(D^{(k)}_{e_{p}}(D^{(k)}_{e_{p}}\psi)-D^{(k)}_{D_{e_{p}}e_{p}}\psi\right)l^{-2}\nonumber
\end{eqnarray}
\end{proof}

Dans tout ce qui suit,  nous noterons $D$ et $\D$ respectivement la dérivée covariante et l'opérateur de Dirac induits par la structure de Weyl $D$ sur le fibré des spineurs de poids conforme $\frac{2-n}{2}$.

\begin{thm}\label{thm2}
Soit $D$ une structure de Weyl. Soient $D$ et $\D$ respectivement la dérivée covariante et l'opérateur de Dirac à poids agissant sur le fibré des spineurs à poids $\Sigma^{\left(\frac{2-n}{2}\right)}$. Pour toutes sections $\psi$ de $\Sigma^{\left(\frac{2-n}{2}\right)}$, nous avons la formule de Lichnerowicz conforme I suivante :
\begin{equation}
\D^{2}\psi=\mathrm{tr}(D^{2})(\psi)+\frac{1}{4}Scal^{D}\psi\nonumber
\end{equation}
\end{thm}

\begin{rem}\label{rem1}
Dans le théorème \ref{thm2}, nous omettons le symbole $\otimes$ du produit tensoriel $Scal^{D}\psi$ de la section $Scal^{D}$ du fibré $L^{(-2)}$ et de la section $\psi$ de $\Sigma^{\left(\frac{2-n}{2}\right)}$. En rappelant que $\Sigma^{(k)}\otimes L^{(-2)}\cong\Sigma^{(k-2)}$, nous remarquons que $Scal^{D}\psi$  est bien une section de $\Sigma^{(k-2)}$. La formule donnée par ce théorème relie bien des objets de même poids conforme $-1-\frac{n}{2}$.
\end{rem}

\begin{proof}
Soient $k\in\N$, $\psi$ et $\phi$ dans $\Sigma^{(k)}$, $x\in M$ et $\{e_{i}\}$ une base de $T_{x}M$. Fixons une métrique $g$ dans $c$. Supposons que la base $\{e_{i}\}$ est parallèle en un point pour la connexion de Levi-Civita $D^{g}$. Nous avons en ce point :
\begin{eqnarray}
D_{e_{i}}^{(k)}(D_{e_{i}}^{(k)}\psi)
&=&D^{g}_{e_{i}}(D^{(k)}_{e_{i}}\psi)+(k-\frac{1}{2})\theta_{g}(e_{i})D_{e_{i}}^{(k)}\psi-\frac{1}{2}e_{i}\cdot\theta_{g}\cdot D_{e_{i}}^{(k)}\psi\nonumber\\
&=&D_{e_{i}}^{g}(D_{e_{i}}^{g}\psi)+(k-\frac{1}{2})D_{e_{i}}^{g}(\theta_{g}(e_{i}))\psi+(2k-1)\theta_{g}(e_{i})D_{e_{i}}^{g}\psi-\frac{1}{2}e_{i}\cdot D_{e_{i}}^{g}(\theta_{g})\cdot\psi\nonumber\\
&&-e_{i}\cdot\theta_{g}\cdot D^{g}_{e_{i}}\psi+(k-\frac{1}{2})^{2}\theta_{g}(e_{i})^{2}\psi
-(k-\frac{1}{2})\theta_{g}(e_{i})e_{i}\cdot\theta_{g}\cdot\psi\nonumber\\
&&+\frac{1}{4}e_{i}\cdot\theta_{g}\cdot e_{i}\cdot\theta_{g}\cdot\psi\nonumber
\end{eqnarray}
Les propriétés de la multiplication de Clifford nous donnent alors :
\begin{eqnarray}
D_{e_{i}}^{(k)}(D_{e_{i}}^{(k)}\psi)
&=&D_{e_{i}}^{g}(D_{e_{i}}^{g}\psi)+(k-\frac{1}{2})D_{e_{i}}^{g}(\theta_{g}(e_{i}))\psi+(2k+1)\theta_{g}(e_{i})D_{e_{i}}^{g}\psi\nonumber\\
&&-\frac{1}{2}e_{i}\cdot D_{e_{i}}^{g}(\theta_{g})\cdot\psi+\theta_{g}\cdot e_{i}\cdot D^{g}_{e_{i}}\psi+(k-\frac{1}{2})^{2}\theta_{g}(e_{i})^{2}\psi\nonumber\\
&&-k\theta_{g}(e_{i})e_{i}\cdot\theta_{g}\cdot\psi-\frac{1}{4}|\theta_{g}|_{g}^{2}\psi\nonumber
\end{eqnarray}
Nous obtenons :
\begin{eqnarray}\label{1}
\sum_{i=1}^{n}D_{e_{i}}^{(k)}(D_{e_{i}}^{(k)}\psi)
&=&-\mathrm{tr}_{g}(D^{g,2})(\psi)
+(k-\frac{1}{2})\mathrm{tr}(D^{g}\theta_{g})\psi
+(2k+1)D_{\theta_{g}^{\sharp}}^{g}\psi-\frac{1}{2}\D^{g}(\theta_{g})\cdot\psi
\nonumber\\
&&+\theta_{g}\cdot\D^{g}\psi+(k^{2}-\frac{n-1}{4})|\theta_{g}|^{2}\psi
\end{eqnarray}
De plus,
$$D_{e_{i}}e_{i}=2\theta_{g}(e_{i})e_{i}-\theta_{g}^{\sharp}$$
En sommant sur $i$ de $1$ à $n$, nous en déduisons :
$$\sum_{i=1}^{n}D_{e_{i}}e_{i}=(2-n)\theta_{g}^{\sharp}$$
Par conséquent,
\begin{equation}\label{2}
\sum_{i=1}^{n}D^{(k)}_{D_{e_{i}}e_{i}}\psi=(2-n)D^{g}_{\theta_{g}^{\sharp}}\psi+k(2-n)|\theta_{g}|_{g}^{2}\psi
\end{equation}
La différence des équations (\ref{1}) et (\ref{2}) nous donne :
\begin{eqnarray}\label{eq5}
-\mathrm{tr}(D^{(k),2})(\psi)&=&-\mathrm{tr}_{g}(D^{g,2})(\psi)+(k-\frac{1}{2})\mathrm{tr}(D^{g}\theta_{g})\psi+(2k+n-1)D^{g}_{\theta_{g}^{\sharp}}\psi-\frac{1}{2}\D^{g}(\theta_{g})\cdot\psi\nonumber\\
&&+\theta_{g}\cdot\D^{g}\psi+\left(k^{2}-k(2-n)-\frac{1}{4}(n-1)\right)|\theta_{g}|_{g}^{2}\psi
\end{eqnarray}
En sommant les expressions (\ref{prop6}) et (\ref{eq5}), puis en utilisant la formule de Lichnerowicz (proposition \ref{prop5}), nous obtenons la formule suivante :
\begin{eqnarray}\label{eq6}
\D^{(k-1)}\D^{(k)}\psi-\mathrm{tr}(D^{(k),2})(\psi)&=&\frac{1}{4}\Big[Scal^{g}+2(2k-1)\mathrm{tr}_{g}(D^{g}(\theta_{g}))-(n-1)(n-2)|\theta_{g}|_{g}^{2}\Big]\psi\nonumber\\
&&+\Big(k+\frac{n-2}{2}\Big)\D^{g}(\theta_{g})\cdot\psi
\end{eqnarray}
Rappelons que la courbure scalaire de $D$, via la métrique $g$, s'identifie à une fonction sur $M$ de la façon suivante :
$$Scal^{D}=Scal^{g}-2(n-1)\mathrm{tr}_{g}(D^{g}\theta_{g})-(n-1)(n-2)|\theta_{g}|_{g}^{2}$$
Ainsi, pour l'unique poids conforme $\frac{2-n}{2}$, la formule (\ref{eq4}) ne dépend pas de la métrique $g$ choisie et nous donne la formule souhaitée :
$$\D^{2}\psi-\mathrm{tr}(D^{2})(\psi)=\frac{1}{4}Scal^{D}\psi$$
\end{proof}

\subsection{Formule de Lichnerowicz conforme II}

Nous allons introduire des opérateurs invariants conformes généralisant la notion de divergence du cas riemannien et calculer des opérateurs adjoints dans un sens conforme. Notons d'abord qu'une structure de Weyl préserve l'application $h$.

\begin{prop}\label{prop8}
Soient $\psi\in\Sigma^{(k)}$ et $\varphi\in\Sigma^{(l)}$. Toute structure de Weyl $D$ préserve la forme bilinéaire $h$ :
\begin{equation}\label{eq7}
\nabla^{D}\left(h(\psi,\varphi)\right)=h(D^{(k)}\psi,\varphi)+h(\psi,D^{(l)}\varphi)
\end{equation}
où $\nabla^{D}$ agit sur $L^{(k+l)}$.
\end{prop}

\begin{proof}
L'application $h$ ne dépend que de la classe conforme $c$. De plus, toute structure de Weyl préserve la classe conforme. Par conséquent, toute structure de Weyl $D$ préserve $h$, et donc nous obtenons la formule (\ref{eq2}).
\end{proof}

Soient $g$ et $\tld{g}$ deux métriques riemanniennes telles que $\tld{g}=f^{2}g$. Notons $\ast^{\tld{g}}$ et $\ast^{g}$ les opérateurs de Hodge associés respectivement aux métriques $\tld{g}$ et $g$. Nous avons : 
$$\ast^{\tld{g}}\omega^{q}=f^{n-2k}\ast^{g}\omega^{q}$$
où $\omega^{q}$ est une $q$-forme sur $M$. Par conséquent, il existe un opérateur de Hodge conforme naturel : 
$$\ast\ :\ \Lambda^{q}M\otimes L^{(k)}\rightarrow\Lambda^{n-q}M\otimes L^{(n-2q+k)}$$
Pour $X$ dans $TM$ et $\omega$ une section de $\Lambda^{q}M\otimes L^{(k)}$, l'opérateur de Hodge conforme vérifie la relation suivante :
\begin{equation}\label{eq8}
X\lrcorner\,\omega=\ast(X^{\flat}\wedge\ast\omega)
\end{equation}
Pour $k=2q-n$, nous pouvons définir un opérateur de divergence conforme par :  
$$\delta=-\ast\circ\,d\circ\ast$$
où $d$ est la différentielle extérieure sur $M$. Nous pouvons également définir un opérateur de divergence relatif à la structure de Weyl $D$, que nous notons $\delta^{D}$. Soit $\{e_{i}\}$ une base locale $c$-orthonormée telle que $c(e_{i},e_{j})=\delta_{ij}l^{2}$. Nous définissons $\delta^{D}$ par :
$$\delta^{D}=-\sum_{i=1}^{n}\left(e_{i}\lrcorner\,D_{e_{i}}\right)l^{-2}$$
Nous pouvons considérer qu'une structure de Weyl $D$ agit sur les sections du fibré $\Lambda^{q}M\otimes L^{(k)}$. En effet, $D$ agit par extension sur l'espace des $q$-formes sur $M$ et sur les fibrés en droites $L^{(k)}$ via sa connexion linéaire associée $\nabla^{D}$. Nous avons alors : 
$$\delta^{D}\ :\ \Lambda^{q}M\otimes L^{(k)}\rightarrow\Lambda^{q-1}M\otimes L^{(k-2)}$$
Lorsqu'une métrique est fixée dans la classe conforme, l'opérateur de Hodge conforme et la divergence conforme s'identifient respectivement avec l'opérateur de Hodge et la divergence riemannienne relative à cette métrique. En revanche, la structure de Weyl n'est pas en général  la connexion de Levi-Civita d'une métrique, par conséquent, les opérateurs $\delta$ et $\delta^{D}$ sont en général distincts.

Notons $D^{a}$ l'opérateur antisymétrisé de la connexion sans torsion $D$ défini par : 
$$D^{a}=\sum_{i=1}^{n}e_{i}^{\ast}\wedge D_{e_{i}}$$
où $\{e_{i}\}$ est une base quelconque de $TM$.

\begin{prop}\label{prop9}
L'opérateur $\delta^{D}$ : $\Lambda^{q}M\otimes L^{(k)}\rightarrow\Lambda^{q-1}M\otimes L^{(k-2)}$ satisfait la relation suivante :
$$\delta^{D}=-\ast\circ D^{a}\circ\ast$$
\end{prop}

\begin{proof}
Soient $\{e_{i}\}$ une base locale $c$-orthonormée telle que $c(e_{i},e_{j})=\delta_{ij}l^{2}$ et $\omega$ une section de $\Lambda^{q}M\otimes L^{(k)}$. D'après la relation (\ref{eq8}), nous avons :
$$\delta^{D}(\omega)=-\sum_{i=1}^{n}e_{i}\lrcorner\,D_{e_{i}}(\omega) l^{-2}
=-\sum_{i=1}^{n}\ast\big(e_{i}^{\flat}\wedge\ast D_{e_{i}}(\omega)\big)l^{-2}$$
De plus, nous avons $e^{\ast}_{i}=e_{i}^{\flat}l^{-2}$, donc :
\begin{eqnarray}
\delta^{D}(\omega)&=&-\sum_{i=1}^{n}\ast\big(e_{i}^{\flat}l^{-2}\wedge\ast D_{e_{i}}(\omega)\big)\nonumber\\
&=&-\sum_{i=1}^{n}\ast\big(e_{i}^{\ast}\wedge\ast D_{e_{i}}(\omega)\big)\nonumber
\end{eqnarray}
Enfin, la structure de Weyl préserve la classe conforme $c$, de telle sorte que l'opérateur $\ast$, qui ne dépend que de $c$, commute avec la connexion $D$. Ainsi :
$$\delta^{D}(\omega)=-\ast\Big(\sum_{i=1}^{n}e_{i}^{\ast}\wedge D_{e_{i}}(\ast\omega)\Big)$$
Nous obtenons la formule suivante :
$$\delta^{D}(\omega)=-\ast\big(D^{a}(\ast\omega)\big)$$
\end{proof}

\begin{cor}\label{prop10}
Sur l'espace des sections du fibré $\Lambda^{q}\otimes L^{(2q-n)}$, l'opérateur de divergence conforme $\delta$ et  l'opérateur de divergence $\delta^{D}$ associé à la structure de Weyl $D$ coïncident.
\end{cor}

\begin{proof}
Si $\omega$ est une section de $\Lambda^{q}\otimes L^{(2q-n)}$, alors $\ast\omega$ est une $(n-q)$-forme sur $M$. Cependant, pour toute connexion $D$ sans torsion sur $TM$, nous avons $d=D^{a}$ sur l'espace des formes sur $M$. Ainsi, d'après la proposition \ref{prop9}, nous avons  l'égalité $\delta=\delta^{D}$ sur $\Lambda^{q}\otimes L^{(2q-n)}$.
\end{proof}

\'Etudions l'existence d'un opérateur adjoint pour la connexion $D^{(k)}$ relativement au produit scalaire $H$ sur les spineurs de poids $k$.  Nous avons $D^{(k)}$ : $\Sigma^{(k)}\rightarrow T^{\ast}M\otimes\Sigma^{(k)}$ ; par conséquent, l'adjoint $D^{(k)\ast}$ de $D^{(k)}$ doit agir sur les sections du fibré $T^{\ast}M\otimes\Sigma^{(k)}$. Soient $\psi$ et $\varphi$ dans $\Sigma^{(k)}$, et $\alpha$ dans $T^{\ast}M$. Nous devons trouver  $D^{(k)\ast}$ tel que l'intégrale du terme suivant ait un sens et soit nulle :
\begin{equation}\label{eq10}
h(D^{(k)}\psi,\alpha\otimes\varphi)-h(\psi,D^{(k)\ast}(\alpha\otimes\varphi))
\end{equation}
Les objets $D^{(k)}\psi$ et $\alpha\otimes\varphi$ sont de poids conforme $k-1$ ; par suite,  $h(D^{(k)}\psi,\alpha\otimes\varphi)$ est de poids $k-2$. Ainsi, pour que les deux termes de l'expression (\ref{eq10}) soient de même poids, $D^{\ast}(\alpha\otimes\varphi)$ doit être de poids conforme $k-2$.  Nous en déduisons que l'intégrale de l'expression (\ref{eq10}) a un sens si et seulement si $k=\frac{2-n}{2}$. Rappelons que nous notons $D$ la structure de Weyl de poids $\frac{2-n}{2}$.

\begin{prop}\label{prop11}
La structure de Weyl $D$ de poids $\frac{2-n}{2}$ possède un adjoint formel, $D^{\ast}$ : $T^{\ast}M\otimes\Sigma^{\left(\frac{2-n}{2}\right)}\rightarrow\Sigma^{\left(\frac{-1-n}{2}\right)}$, défini par : 
$$D^{\ast}(\alpha\otimes\varphi)=-D_{\alpha^{\sharp}}\varphi+\delta^{D}(\alpha)\varphi$$
pour $\alpha\in T^{\ast}M$ et $\varphi\in\Sigma^{\left(\frac{2-n}{2}\right)}$. De plus, si $\psi$ est une section de $\Sigma^{\left(\frac{2-n}{2}\right)}$, nous avons la formule suivante : 
\begin{equation}\label{eqprop11}
h\big(D\psi,\alpha\otimes\varphi\big)-h\big(\psi,D^{\ast}(\alpha\otimes\varphi)\big)=-\delta^{D}\big(h(\psi,\alpha\otimes\varphi)\big)
\end{equation}
Cette formule relie des sections du fibré $L_{\C}^{(-n)}$
\end{prop}

\begin{proof}
Soient $\psi$ et $\varphi$ deux sections de $\Sigma^{\left(\frac{2-n}{2}\right)}$ à support compact, et $\alpha$ une $1$-forme sur $M$. Calculons 
$$H(D\psi,\alpha\otimes\varphi)=\int_{M}h(D\psi,\alpha\otimes\varphi)$$
Soit $\{e_{i}\}$ une base locale $c$-orthonormée telle que $c(e_{i},e_{j})=\delta_{ij}l^{2}$. Nous avons alors : 
\begin{eqnarray}
h(D_{e_{i}}\psi,\alpha(e_{i})\varphi)l^{-2}
&=&\nabla_{e_{i}}^{D}\left(h(\psi,\alpha(e_{i})\varphi) l^{-2}\right)-h\big(\psi,D_{e_{i}}(\alpha(e_{i})\varphi l^{-2})\big)\nonumber\\
&=&\nabla_{e_{i}}^{D}\big(h(\psi,\alpha(e_{i})\varphi)\big)l^{-2}-h(\psi,\alpha(D_{e_{i}}e_{i})\varphi)l^{-2}\nonumber\\
&&-h\big(\psi,D_{e_{i}}(\alpha)(e_{i})l^{-2}\varphi\big)-h\big(\psi,\alpha(e_{i})l^{-2}D_{e_{i}}\varphi\big)\nonumber
\end{eqnarray}
En sommant sur $i$, nous obtenons :
\begin{eqnarray}
\sum_{i=1}^{n}h(D_{e_{i}}\psi,\alpha(e_{i})\varphi)l^{-2}&=&\sum_{i=1}^{n}\Big(\nabla_{e_{i}}^{D}\big(h(\psi,\alpha(e_{i})\varphi)\big)-h\big(\psi,\alpha(D_{e_{i}}e_{i})\varphi\big)\Big)l^{-2}\nonumber\\
&&+h\Big(\psi,-\sum_{i=1}^{n}\left(e_{i}\lrcorner\,D_{e_{i}}(\alpha)l^{-2}\right)\varphi\Big)-h\Big(\psi,\sum_{i=1}^{n}\alpha(e_{i})D_{e_{i}}\varphi\Big)\nonumber
\end{eqnarray}
Nous avons donc la formule suivante :
$$h\big(D\psi,\alpha\otimes\varphi\big)-h\big(\psi,(-D_{\alpha^{\sharp}}+\delta^{D}(\alpha))\varphi\big)=-\delta^{D}\big(h(\psi,\alpha\otimes\varphi)\big)$$
Cependant, $h(\psi,\alpha\otimes\varphi)$ étant une section de $\Lambda^{1}M\otimes L^{(2-n)}$, le corollaire \ref{prop10} montre que :
$$\delta^{D}\big(h(\psi,\alpha\otimes\varphi)\big)=\delta\big(h(\psi,\alpha\otimes\varphi)\big)$$
Ainsi, d'après la formule de Stokes, nous avons : 
$$\int_{M}\delta^{D}\big(h(\psi,\alpha\otimes\varphi)\big)=0$$
Nous obtenons alors : 
$$H(D\psi,\alpha\otimes\varphi)=H\big(\psi,[-D_{\alpha^{\sharp}}+\delta^{D}(\alpha)]\varphi\big)$$
L'adjoint de $D$ est donc donné par la formule suivante : 
$$D^{\ast}(\alpha\otimes\psi)=-D_{\alpha^{\sharp}}\psi+\delta^{D}(\alpha)\psi$$
\end{proof}

\begin{lem}\label{lem1}
Soient $\{e_{i}\}$ une base $c$-orthonormée telle que $c(e_{i},e_{j})=\delta_{ij}l^{2}$. Nous avons la formule suivante : 
$$\sum_{i=1}^{n}\delta^{D}(e_{i}^{\ast})D_{e_{i}}\psi=\sum_{i=1}^{n}l^{-2}D_{D_{e_{i}}e_{i}}\psi$$
pour tout $\psi$ dans $\Sigma^{\left(\frac{2-n}{2}\right)}$.
\end{lem}

\begin{proof}
Dans un premier temps, nous calculons : 
\begin{eqnarray}
\delta^{D}(e_{i}^{\ast})
&=&-\sum_{j=1}^{n}e_{j}\lrcorner\,D_{e_{j}}(e_{i}^{\ast})l^{-2}\nonumber\\
&=&-\sum_{j=1}^{n}D_{e_{j}}(e_{i}^{\ast}(e_{j}))l^{-2}+\sum_{j=1}^{n}e_{i}^{\ast}\big(D_{e_{j}}e_{j}\big)l^{-2}\nonumber\\
&=&-\sum_{j=1}^{n}D_{e_{j}}(\delta_{ij})l^{-2}+\sum_{j=1}^{n}e_{i}^{\ast}\big(D_{e_{j}}e_{j}\big)l^{-2}\nonumber\\
&=&\sum_{j=1}^{n}e_{i}^{\ast}\big(D_{e_{j}}e_{j}\big)l^{-2}\nonumber
\end{eqnarray}
D'autre part, nous avons : 
$$\sum_{i=1}^{n}\delta^{D}(e_{i}^{\ast})D_{e_{i}}\psi=\sum_{i=1}^{n}D_{\delta^{D}(e_{i}^{\ast})e_{i}}\psi$$
D'après ce qui précède, nous obtenons : 
$$\sum_{i=1}^{n}\delta^{D}(e_{i}^{\ast})D_{e_{i}}\psi=\sum_{j=1}^{n}l^{-2}D_{\sum_{i=1}^{n}e_{i}^{\ast}\big(D_{e_{j}}e_{j}\big)e_{i}}\psi$$
Enfin, nous avons $\sum_{i=1}^{n}e_{i}^{\ast}\big(D_{e_{j}}e_{j}\big)e_{i}=D_{e_{j}}e_{j}$, ce qui nous donne :
$$\sum_{i=1}^{n}\delta^{D}(e_{i}^{\ast})D_{e_{i}}\psi=\sum_{i=1}^{n}l^{-2}D_{D_{e_{i}}e_{i}}\psi$$
\end{proof}

\begin{prop}\label{prop12}
Soit $\psi$ un spineur de poids $\frac{2-n}{2}$. Pour toute structure de Weyl $D$, nous avons la formule suivante :
$$D^{\ast}D\psi=\mathrm{tr}(D^{2})(\psi)$$
où $D$ agit sur $\Sigma^{\left(\frac{2-n}{2}\right)}$ via la structure de Weyl de poids correspondant.
\end{prop}

\begin{proof}
Par définition de l'adjoint formel de la connexion $D$, nous avons :
$$D^{\ast}D\psi=D^{\ast}\Big(\sum_{i=1}^{n}e_{i}^{\ast}\otimes D_{e_{i}}\psi\Big)
=-\sum_{i=1}^{n}\left(D_{(e_{i}^{\ast})^{\sharp}}(D_{e_{i}}\psi)-\delta^{D}(e_{i}^{\ast})D_{e_{i}}\psi\right)$$
Cependant, nous avons $(e_{i}^{\ast})^{\sharp}=e_{i}l^{-2}$, et donc, d'après le lemme \ref{lem1}, nous obtenons :
$$D^{\ast}D\psi=-\sum_{i=1}^{n}\big(D_{e_{i}}(D_{e_{i}}\psi)-D_{D_{e_{i}}e_{i}}\psi\big)l^{-2}$$
\end{proof}

Nous allons maintenant nous intéresser à des formule de type Lichnerowicz conforme mettant en jeux des sections du fibré des densités $L^{(-n)}$.

\begin{prop}\label{prop13}
Soient $\psi$ et $\varphi$ deux sections du fibré $\Sigma^{\left(\frac{2-n}{2}\right)}$. Nous avons la formule suivante :
$$h(D\psi,D\varphi)=h(\psi,\mathrm{tr}(D^{2})(\varphi))-\delta^{D}\big(h(\psi,D\varphi)\big)\in L_{\C}^{(-n)}$$
\end{prop}

\begin{proof}
Remarquons que $D\varphi$ est une section de $T^{\ast}M\otimes\Sigma^{\left(\frac{2-n}{2}\right)}$. En remplaçant $\alpha\otimes\varphi$ par $D\varphi$ dans la formule (\ref{eqprop11}) de la proposition \ref{prop11}, nous obtenons :
$$h(D\psi,D\varphi)-h(\psi,D^{\ast}D\varphi)=-\delta^{D}\big(h(\psi,D\varphi)\big)$$
La proposition \ref{prop12} permet de conclure :
$$h\big(D\psi,D\varphi\big)-h\big(\psi,\mathrm{tr}(D^{2})(\varphi)\big)=-\delta^{D}\big(h(\psi,D\varphi)\big)$$
\end{proof}

\begin{cor}\label{cor13}
Soient $\psi$ et $\varphi$ deux sections du fibré $\Sigma^{\left(\frac{2-n}{2}\right)}$. Nous avons la formule suivante :
$$h(D\psi,D\varphi)+\frac{1}{4}h(\psi,Scal^{D}\varphi)-h(\psi,\D^{2}\varphi)=-\delta^{D}\big(h(\psi,D\varphi)\big)\in L_{\C}^{(-n)}$$
\end{cor}

\begin{proof}
La proposition \ref{prop13} nous donne : 
$$h(D\psi,D\varphi)+\frac{1}{4}h(\psi,Scal^{D}\varphi)-h(\psi,\D^{2}\varphi)=h(\psi,\mathrm{tr}(D^{2})(\varphi)+\frac{1}{4}Scal^{D}\varphi-\D^{2}\varphi)-\delta^{D}\big(h(\psi,D\varphi)\big)$$
De plus, la formule de Lichnerowicz conforme I (théorème \ref{thm2}) nous donne : 
$$\mathrm{tr}(D^{2})(\varphi)+\frac{1}{4}Scal^{D}\varphi-\D^{2}\varphi=0$$
Nous obtenons donc la formule souhaitée.
\end{proof}

Menons le même calcul que dans la proposition \ref{prop13} pour l'opérateur de Dirac $\D$ de poids conforme $\frac{2-n}{2}$. 

\begin{prop}\label{prop14}
Soient $\psi$ et $\varphi$ deux sections du fibré $\Sigma^{\left(\frac{2-n}{2}\right)}$. Nous avons la formule suivante :
$$h(\D\psi,\D\varphi)=h(\psi,\D^{2}\varphi)+\delta^{D}(\beta_{\psi,\varphi})\in L_{\C}^{(-n)}$$
où $\beta_{\psi,\varphi}$ est la section de $T^{\ast}M\otimes L^{(2-n)}$ définie par : $\beta_{\psi,\varphi}(X)=h(\psi,X^{\flat}\cdot\D\varphi)$
\end{prop}

\begin{proof}
Soit $\{e_{i}\}$ une base $c$-orthonormée telle que $c(e_{i},e_{j})=\delta_{ij}l^{2}$. Les propriétés de l'application $h$ et de la structure de Weyl $D$ nous donnent : 
\begin{eqnarray}
h(\D\psi,\D\varphi)&=&\sum_{i=1}^{n}h\big(e_{i}^{\ast}\cdot D_{e_{i}}\psi,\D\varphi\big)\nonumber\\
&=&-\sum_{i=1}^{n}h\big(D_{e_{i}}\psi,e_{i}^{\ast}\cdot\D\varphi\big)\nonumber\\
&=&-\sum_{i=1}^{n}\Big(\nabla^{D}\big(h(\psi,e_{i}^{\ast}\cdot\D\varphi)\big)-h\big(\psi,D_{e_{i}}(e_{i}^{\ast}\cdot\D\varphi)\big)\Big)\nonumber\\
&=&-\sum_{i=1}^{n}\Big(\nabla^{D}\big(h(\psi,e_{i}^{\ast}\cdot\D\varphi)\big)-h\big(\psi,D_{e_{i}}(e_{i}^{\ast})\cdot\D\varphi\big)\Big)+\sum_{i=1}^{n}h\big(\psi,e_{i}^{\ast}\cdot D_{e_{i}}(\D\psi)\big)\nonumber\\
&=&-\sum_{i=1}^{n}\Big(\nabla^{D}\big(h(\psi,e_{i}^{\flat}l^{-2}\cdot\D\varphi)\big)-h\big(\psi,D_{e_{i}}(e_{i}^{\flat}l^{-2})\cdot\D\varphi\big)\Big)+h(\psi,\D^{2}\varphi)\nonumber\\
&=&-\sum_{i=1}^{n}\Big(\nabla^{D}\big(h(\psi,e_{i}^{\flat}\cdot\D\varphi)\big)-h\big(\psi,D_{e_{i}}(e_{i}^{\flat})\cdot\D\varphi\big)\Big)l^{-2}+h(\psi,\D^{2}\varphi)\nonumber\\
&=&-\sum_{i=1}^{n}\Big(\nabla^{D}\big(h(\psi,e_{i}^{\flat}\cdot\D\varphi)\big)-h\big(\psi,(D_{e_{i}}e_{i})^{\flat}\cdot\D\varphi\big)\Big)l^{-2}+h(\psi,\D^{2}\varphi)\nonumber
\end{eqnarray}
En posant $\beta_{\psi,\varphi}(X)=h(\psi,X^{\flat}\cdot\D\varphi)$, nous obtenons la formule souhaitée :
$$h(\D\psi,\D\varphi)=h(\psi,\D^{2}\varphi)+\delta^{D}(\beta_{\psi,\varphi})$$
\end{proof}
 
\begin{thm}\label{thm3}
Soit $(M,c)$ une variété conforme orientée. Pour toute structure de Weyl $D$ sur $M$, nous avons la formule de Lichnerowicz conforme II suivante :
$$h(D\psi,D\varphi)+\frac{1}{4}h(\psi,Scal^{D}\varphi)-h(\D\psi,\D\varphi)=-\delta(\omega_{\psi,\varphi})$$
où $\psi$ et $\varphi$ sont des sections de $\Sigma^{\left(\frac{2-n}{2}\right)}$. La section $\omega_{\psi,\varphi}$ du fibré $T^{\ast}M\otimes L^{(2-n)}$ est  définie par : $\omega_{\psi,\varphi}(X)=h\big(\psi,X^{\flat}\cdot\D\varphi+D_{X}\varphi\big)$
\end{thm}

Notons  $|\psi|_{h}^{2}=h\big(\psi,\psi\big)$, pour $\psi\in\Sigma$. D'après le théorème \ref{thm3}, pour toute structure de Weyl $D$ sur la variété conforme $(M,c)$ et toute section $\psi$ du fibré des spineurs à poids $\Sigma$, nous avons :
\begin{equation}\label{eq11}
|\psi|_{h}^{2}+\frac{1}{4}Scal^{D}|\psi|_{h}^{2}-|\D\psi|_{h}^{2}=-\delta(\omega_{\psi})
\end{equation}
où $\omega_{\psi,\psi}$ est noté $\omega_{\psi}$. 

\begin{rem}\label{rem2}
Les objets $D\psi$, $D\varphi$, $\D\psi$ et $\D\varphi$ sont de poids conforme $-n/2$. Par conséquent, $h(D\psi,D\varphi)$ et $h(\D\psi,\D\varphi)$ sont de poids $-n$. De plus, le produit tensoriel $Scal^{D}\varphi$ est de poids conforme $-1-n/2$, donc $h(\psi,Scal^{D}\varphi)$ est aussi de poids $-n$. Enfin, par définition de l'opérateur $\delta$, l'objet $\delta(\omega_{\psi,\varphi})$ est une section de $L_{\C}^{(-n)}$. Le théorème ci-dessus donne une formule reliant des sections du fibré des densités d'intégration $L_{\C}^{(-n)}$.
\end{rem} 

\begin{proof}
Les propositions \ref{prop13} et \ref{prop14} nous donnent :
$$h(D\psi,D\varphi)+\frac{1}{4}h(\psi,Scal^{D}\varphi)-h(\D\psi,\D\varphi)=h\big(\psi,\mathrm{tr}(D^{2})(\varphi)+\frac{1}{4}Scal^{D}\varphi-\D^{2}\varphi\big)-\delta^{D}(\omega_{\psi,\varphi})$$
où $\omega_{\psi,\varphi}=h(\psi,D\varphi)+\beta_{\psi,\varphi}$. D'après la formule de Lichnerowicz conforme I (théorème \ref{thm2}), nous avons :
$$\mathrm{tr}(D^{2})(\varphi)+\frac{1}{4}Scal^{D}\varphi-\D^{2}\varphi=0$$
Nous avons donc la formule suivante : 
$$h(D\psi,D\varphi)+\frac{1}{4}h(\psi,Scal^{D}\varphi)-h(\D\psi,\D\varphi)=-\delta^{D}(\omega_{\psi,\varphi})$$
Le corollaire \ref{prop10} assure l'égalité des opérateurs $\delta^{D}$ et $\delta$ sur les sections du fibré $L^{(2-n)}$. Ainsi, $\omega_{\psi,\varphi}$ étant une section de $L^{(2-n)}$, nous obtenons la formule souhaitée. 
\end{proof}

\section{Structures conformes asymptotiquement plates}

Soit $(M,g)$ une variété riemannienne orientée, complète et non compacte de dimension $n$. Notons $g_{can}$ la métrique plate sur $\R^{n}$.

\subsection{Variétés asymptotiquement plates}

Notons $E_{R}=\R^{n}\setminus B_{R}$ l'extérieur de la boule de centre $0$ et de rayon $R$ de $\R^{n}$.
Supposons qu'il existe un compact $K$ de $M$ et un difféomorphisme
$$\Phi\,:\,E_{R}\rightarrow M\setminus K$$ Notons $V=M\setminus K$. Le couple $(V,\Phi)$ est une carte
à l'infini et $V$ le bout de $M$. Dans cette carte $\Phi(x)=(x_{1},\ldots,x_{n})$, la norme d'un élément s'écrit :
$$|x|=r=\left(\sum_{i=1}^{n}x_{i}^{2}\right)^{\frac{1}{2}}$$

Soit $r>R$; notons $M_{r}=M\setminus E_{r}$, où $E_{r}$ est confondu avec son image par $\Phi$. L'ensemble $M_{r}$ est un compact de $M$ dont le bord s'identifie à la sphère $S_{r}$ de $\R^{n}$.

\begin{defn}\label{defn3}
Une variété riemannienne $(M,g)$ est asymptotiquement plate d'ordre $\tau>0$ s'il existe une décomposition $M\setminus K=\sqcup_{l=1}^{k} M_{\infty}^{l}$, où $K$ est un compact, et des difféomorphismes $\Phi_{l}$ de $M_{\infty}^{l}$ dans $E_{R_{l}}$, pour $R_{l}>0$, tels que :
$$g_{ij}=\delta_{ij}+O(r_{l}^{-\tau})\,,\ \partial_{k}g_{ij}=O(r_{l}^{-\tau-1})\,et\ \partial_{l}\partial_{k}g_{ij}=O(r_{l}^{-\tau-2})\,,$$
quand $r_{l}=|x_{l}|\rightarrow\infty$, dans les coordonnées engendrées par le difféomorphisme $\Phi_{l}$ sur $M_{\infty}^{l}$, pour $l$ de $1$ à $k$. Les ouverts $M_{\infty}^{l}$ sont les bouts de $M$ et les coordonnées $\{x_{l}^{i}\}$ sont appelées coordonnées asymptotiques sur $M_{\infty}^{l}$.
\end{defn}

Supposons que $(M,g)$ est une variété asymptotiquement plate d'ordre $\tau$  ne possédant qu'un seul bout $M_{\infty}$. Soient $p>1$ et $\delta\in\R$. Nous notons $\nabla$ la connexion de Levi-Civita associée à $g$. Sauf mention explicite du contraire, tout les objets riemanniens considérés sont relatifs à la métrique $g$.
L'espace des fonctions $L^{p}_{\delta}$ est défini comme le complété de l'espace des fonctions $C^{\infty}$ à support compact sur $M$ pour la norme :
$$||u||_{p,\delta}=\left(\int_{M}|u|^{p}r^{-\delta p-n}v_{g}\right)^{1/p}$$
où $v_{g}$ est la forme volume associée à $g$ et $r$ est la distance radiale sur le bout $M_{\infty}$ prolongée par $1$ sur la partie compacte de $M$.
Nous définissons également le $k$-ième espace de Sobolev de poids $\delta$ sur $M$, que l'on note $W^{k,p}_{\delta}$, comme le complété de l'espace des fonctions $C^{\infty}$ à support compact sur $M$ pour la norme de Sobolev suivante :
\begin{equation}\label{eqsobo}
||u||_{k,p,\delta}=\sum_{j=0}^{k}||\nabla^{j}u||_{p,\delta-j}
\end{equation}
L'espace des fonctions $C_{\delta}^{k}$ est défini comme l'ensemble des fonctions $u$, $C^{k}$ sur $M$, dont la norme suivante est finie :
\begin{equation}\label{eqholder1}
||u||_{C^{k}_{\delta}}=\sum_{j=0}^{k}\sup_{M_{\infty}}r^{-\delta+j}|\nabla^{j}u|
\end{equation}
Soit $\alpha\in]0,1[$; l'espace de Hölder de poids $\delta$ est défini par l'ensemble des fonctions $u$ dans $C^{k}_{\delta}$ telles que la norme suivante est finie :
\begin{equation}\label{eqholder2}
||u||_{C^{k,\alpha}_{\delta}}=||u||_{C^{k}_{\delta}}+\sup_{x,y\in M_{\infty}}\left(\min(r(x),r(y))^{-\delta+k+\alpha}
\frac{|\nabla^{k}u(x)-\nabla^{k}u(y)|}{|x-y|^{\alpha}}\right)
\end{equation}
où $y$ est dans un voisinage de $x$ et $\nabla^{k}u(y)$ est le tenseur en $x$ obtenu par transport parallèle le long de la géodésique radiale joignant $x$ à $y$. Ces espaces de fonctions dépendent des coordonnées choisies sur le bout $M_{\infty}$ de la variété. En revanche, les différents systèmes de coordonnées étant asymptotiques aux coordonnées euclidiennes, les normes relatives à deux systèmes de coordonnées différents sont équivalentes. La définition des normes de Hölder (\ref{eqholder1}) et (\ref{eqholder2}) nous donne immédiatement la proposition suivante :

\begin{prop}\label{prop15}
 Soit $u$ appartenant à $C^{2,\alpha}_{-\tau}(M)$. Alors $u=O(r^{-\tau})$ et, pour $i$ et $j$ dans $\{1,\ldots,n\}$, $\partial_{i}u=O(r^{-\tau-1})$ et $\partial^{2}_{ij}u=O(r^{-\tau-2})$.
\end{prop}

Pour ces espaces de fonctions, nous avons également un théorème de Sobolev à poids \cite{lp} :

\begin{prop}\label{prop16}
Soient $q>1$ et $\alpha\in]0,1[$. Supposons que $l-k-\alpha>\frac{n}{q}$. Pour tout $\varepsilon>0$, nous avons les inclusions continues suivantes :
$$C_{\delta-\varepsilon}^{l,\alpha}\subset W_{\delta}^{l,q}\subset C_{\delta}^{k,\alpha}$$
En particulier, si $u\in W_{\delta}^{l,q}$ avec $l>\frac{n}{q}$, $u=O(r^{\delta})$.
\end{prop}

\begin{defn}\label{defn4}
La variété riemannienne $(M,g)$ asymptotiquement plate d'ordre $\tau$ vérifie les conditions de décroissance de la masse lorsque :
\begin{itemize}
\item l'ordre $\tau$ de $M$ est tel que :
$$\frac{n-2}{2}<\tau<n-2$$
\item le tenseur $g-g_{can}$ appartient à $C_{-\tau}^{2,\alpha}(M_{\infty}^{l})$, pour tout $l$ de $1$ à $k$.
\item la courbure scalaire $Scal^{g}$ de la connexion de Levi-Civita $\nabla$ est intégrable sur $M$.
\end{itemize}
Notons $\M_{\tau}$ l'espace des métriques sur $M$ vérifiant les conditions ci-dessus.
\end{defn}

Lorsque une variété asymptotiquement plate $(M,g)$ d'ordre $\tau$ à un seul bout $M_{\infty}$ vérifie les conditions de la définition précédente, l'expression de la masse de cette variété est donnée par :
$$m(g)=\lim_{r\rightarrow\infty}\int_{S_{r}}\sum_{i,j=1}^{n}(\partial_{i}g_{ij}-\partial_{j}g_{ii})e_{j}\lrcorner v_{g}$$
où $\nu$ est le champ de vecteurs unitaire sortant de la sphère $S_{r}$ et $g_{ij}$ sont les composantes de la métrique $g$ dans les coordonnées asymptotiques sur $M_{\infty}$.
Si $g$ appartient à l'espace $\M_{\tau}$, la masse de la variété riemannienne $(M,g)$ est bien définie et ne dépend pas du système de coordonnées choisi \cite{bar}. Nous allons rappeler sans démonstration les résultats analytiques nécessaires à la démonstration du théorème de la masse positive dans le cas d'une variété spinorielle \cite{bar}. Ces résultats s'appuient sur les propriétés des espaces de Sobolev, de Hölder et sur des résultats de théorie elliptique pour lesquels le lecteur pourra se référer à \cite{bar}. Fixons une métrique $g$ dans $\M_{\tau}$ et supposons que la variété $M$ est spinorielle. Soit $q$ un entier strictement supérieur à $n$.

\begin{prop}\label{prop17}\emph{(\cite{bar} page 676)}
Soient $p$ dans $]1,+\infty[$ et $\delta$ non exceptionnel \footnote{Le paramètre $\delta\in\R$ est dit non exceptionnel s'il appartient à l'ensemble $\R\setminus\{k\in\Z\,:\,k\neq-1,-2,\cdots,3-n\}$.}. L'opérateur Laplacien
$\Delta$ : $W^{2,p}_{\delta}\rightarrow L^{p}_{\delta-1}$ est un opérateur de Fredholm. De plus, si $2-n<\delta<0$, c'est un isomorphisme de $W^{2,p}_{\delta}$ sur $L^{p}_{\delta-1}$.
\end{prop}

Nous définissons $W^{k,q}_{\delta}(\Sigma^{g})$ l'espace de Sobolev des sections de $\Sigma^{g}$ dont la norme de Sobolev définie de façon analogue à (\ref{eqsobo}) est finie; lorsqu'il n'y a pas d'ambigüités nous noterons également $W_{\delta}^{k,q}$ ces espaces. Le bout $M_{\infty}$ de $M$, identifié à $\R^{n}\setminus B_{R}$, est munit des coordonnées asymptotiques $x=(x_{1},\ldots,x_{n})$. Soit $e=(e_{1},\ldots,e_{n})$ la base orthonormée de $\R^{n}$ induite par $x$. Nous avons $g(\cdot,\cdot)=g_{can}(A\cdot,\cdot)$, où $A$ est un champ de matrices symétriques définies positives. Nous définissons le repère $g$-orthonormé $s$ de $M_{\infty}$ par : $s=(A^{\frac{1}{2}})^{-1}e$, où $A^{\frac{1}{2}}$ est l'unique racine carré définie positive de $A$. Soit $\tilde{s}$ l'un des deux repères spinoriels relevant $s$. Un spineur $\psi_{0}$ de $\Sigma^{g}$ s'appelle un spineur constant si $\psi_{0}=[\tilde{s},\xi_{0}]$, où la fonction $\xi$ : $M_{\infty}\rightarrow\Delta_{n}$ est constante. En particulier, si $\psi_{0}$ est constant, sa norme $|\psi_{0}|$ est constante sur le bout à l'infini $M_{\infty}$.

\begin{defn}\label{defn5}
Une section $\psi$ de $\Sigma^{g}$ est asymptotiquement constante s'il existe un spineur constant $\psi_{0}$ tel que $\psi-\psi_{0}\in W^{2,q}_{-\tau}(\Sigma^{g})$.
\end{defn}

\begin{prop}\label{prop18}\emph{(\cite{bar} page 690)}
Soit $1-n<\delta<0$. L'opérateur de Dirac $\D^{g}$ : $W^{2,q}_{\delta}(\Sigma^{g})\rightarrow W^{1,q}_{\delta-1}(\Sigma^{g})$ est un isomorphisme.
\end{prop}

Nous pouvons déduire de cette proposition l'existence d'un spineur $\D$-harmonique asymptotiquement constant. L'existence d'un tel spineur joue un rôle majeur dans la preuve de théorème de la masse positive.

\begin{cor}\label{cor19}\emph{(\cite{bar} page 690)}
Soit $\psi_{0}$ un spineur constant sur $M$. Il existe un spineur $\psi$ dans $\Sigma^{g}$ tel que :
$$\D^{g}\psi=0$$
et
$$\psi-\psi_{0}\in W_{-\tau}^{2,q}(\Sigma^{g})$$
\end{cor}

\begin{proof}
Les hypothèses de décroissance asymptotique donne $\D^{g}\psi_{0}\in W_{-\tau-1}^{1,q}(\Sigma^{g})$. D'après la proposition \ref{prop18}, il existe un unique spineur $\psi_{1}\in W_{-\tau}^{2,q}(\Sigma^{g})$ tel que $\D^{g}\psi_{1}=-\D^{g}\psi_{0}$. Ainsi, le spineur $\psi=\psi_{1}+\psi_{0}$ convient.
\end{proof}

Dans la méthode de Witten, en intégrant la formule de Lichnerowicz sur la variété asymptotiquement plate $M$, le terme de divergence converge sur le bout à l'infini vers la masse de la variété. Ce fait est donné par la proposition suivante :

\begin{prop}\label{prop19}\emph{(\cite{bar} page 691)}
Soit $\psi$ dans $\Sigma^{g}$ asymptotiquement constant. Nous avons alors la formule suivante :
$$\int_{M}|\nabla\psi|^{2}v_{g}+\frac{1}{4}\int_{M}Scal^{g}|\psi|^{2}v_{g}-\int_{M}|\D^{g}\psi|^{2}v_{g}=\lim_{r\rightarrow\infty}\int_{S_{r}}(\nabla_{\nu}\psi+\nu\cdot\D^{g}\psi,\psi)_{g}\nu\lrcorner\,v_{g}=\frac{1}{4}m(g)|\psi_{0}|^{2}$$
où $\psi_{0}$ est le spineur constant vérifiant $\psi-\psi_{0}\in W^{2,q}_{-\tau}(\Sigma^{g})$.
\end{prop}

Nous avons alors le théorème de la masse positive \cite{bar} :

\begin{thm}\label{thm4}
Soit $(M,g)$ une variété riemannienne, spinorielle et asymptotiquement plate d'ordre $\tau$ telle que $g\in\M_{\tau}$. Supposons que la courbure scalaire de la connexion de Levi-Civita de $g$ est positive. Alors la masse de $(M,g)$ est positive. De plus, la masse est nulle si et seulement si $M$ est isométrique à l'espace $\R^{n}$ euclidien.
\end{thm}

\begin{proof}
Soit $\psi_{0}$ un spineur constant. La proposition \ref{prop19}, appliquée au spineur donné par le corollaire \ref{cor19}, donne la formule suivante :
$$\int_{M}|\nabla\psi|^{2}v_{g}+\frac{1}{4}\int_{M}Scal^{g}|\psi|^{2}v_{g}=\frac{1}{4}m(g)|\psi_{0}|_{g}^{2}$$
Par hypothèse, la courbure scalaire $Scal^{g}$ de $g$ est positive, donc la masse $m(g)$ de $M$ est positive.
\end{proof}

Considérons l'ensemble de fonctions suivant :
$$\F=\{f\in C^{\infty}(M,\R^{>0})\,:\,f-1\in C^{2,\alpha}_{-\tau}\ et\ \Delta(f)\in L^{1}\}$$
Chaque fonction dans $\F$ nous donne une métrique asymptotiquement plate dans la classe conforme de $g$. Nous avons le théorème suivant \cite{sim} :

\begin{thm}\label{thm5}
Soit $(M,g)$ une variété riemannienne asymptotiquement plate pour laquelle $g\in\M_{\tau}$. Soit $\tilde{g}=fg$. Alors, $\tld{g}$ appartient à $\M_{\tau}$ si et seulement si $f$ appartient à $\F$. De plus, leurs masses sont reliées par la formule suivante :
\begin{equation}\label{mass}
m(\widetilde{g})=m(g)+(n-1)\int_{M}\Delta^{g}(f)v_{g}
\end{equation}
\end{thm}

Nous démontrons uniquement la formule (\ref{mass}). Dans la carte à l'infini associée à $g$, nous avons :
$$\partial_{i}(fg_{ij})-\partial_{j}(fg_{ii})=f(\partial_{i}(g_{ij})-\partial_{j}(g_{ii}))+(g_{ij}\partial_{i}f-g_{ii}\partial_{j}f)$$
Posons $\mu_{j}=\sum_{i=1}^{n}(\partial_{i}(g_{ij})-\partial_{j}(g_{ii}))\nu_{j}$. En écrivant $g_{ij}=\delta_{ij}+a_{ij}$ nous obtenons :
$$\sum_{i=1}^{n}\left(\partial_{i}(fg_{ij})-\partial_{j}(fg_{ii})\right)\nu_{j}=\mu_{j}+(f-1)\mu_{j}+(1-n)\partial_{j}f\nu_{j}+\sum_{i=1}^{n}(a_{ij}\partial_{i}f-a_{ii}\partial_{j}f)\nu_{j}$$
Pour terminer notre calcul, nous devons intégrer sur $S_{r}$ et faire tendre $r$ vers l'infini. Cependant, $\partial_{i}f=O(r^{-\tau-1})$, $a_{ij}=O(r^{-\tau})$, $\mu_{j}=O(r^{-\tau-1})$ et $f-1=O(r^{-\tau})$. Par conséquent
$$\sum_{i=1}^{n}(a_{ij}\partial_{i}f-a_{ii}\partial_{j}f)\nu_{j}=O(r^{-2\tau-1})\,\ \text{et}\ \, (f-1)\mu_{j}=O(r^{-2\tau-1})$$
avec $2\tau+1>n-1$, donc les intégrales sur $S_{r}$ de ces deux termes tendent vers $0$ lorsque $r$ tend vers l'infini. En sommant sur $j$, nous en déduisons la formule suivante :
$$m(\widetilde{g})=m(g)+(1-n)\lim_{r\rightarrow\infty}\int_{S_{r}}df(\nu)\nu\lrcorner v_{g}$$
La formule de Stokes termine la démonstration :
$$m(\widetilde{g})=m(g)+(n-1)\lim_{r\rightarrow\infty}\int_{M_{r}}\Delta_{g}(f)v_{g}$$

Terminons ce paragraphe par un lemme dont nous aurons besoin ultérieurement dans cette note.

\begin{lem}\label{lem2}
Soit $(M,g)$ une variété asymptotiquement plate telle que $g\in\M_{\tau}$.
Pour toute fonction $f$ dans $\F$, nous avons :
$$\lim_{r\rightarrow\infty}\int_{S_{r}}\frac{df}{f}(\nu)\nu\lrcorner\,v_{g}=\lim_{r\rightarrow\infty}\int_{S_{r}}df(\nu)\nu\lrcorner\,v_{g}$$
\end{lem}

\begin{proof}
Nous écrivons :
\begin{equation}\label{eq14}
\int_{S_{r}}\frac{df}{f}(\nu)\nu\lrcorner\,v_{g}=\int_{S_{r}}(f^{-1}-1)df(\nu)\nu\lrcorner\,v_{g}+\int_{S_{r}}df(\nu)\nu\lrcorner\,v_{g}
\end{equation}
La fonction $f$ est strictement positive sur $M$ et tend vers $1$ à l'infini, il existe donc $\varepsilon>0$ et $C>0$ tels que $\varepsilon\leq f\leq C$. Par conséquent, la fonction $f^{-1}$ est bornée sur $M$. De plus, $f-1=O(r^{-\tau})$, nous avons alors : 
$f^{-1}-1=-f^{-1}(f-1)=O(r^{-\tau})$, donc $(f^{-1}-1)|df|=O(r^{-2\tau-1})$ avec $2\tau+1>n-1$. Nous en déduisons : 
$$\lim_{r\rightarrow\infty}\int_{S_{r}}(f^{-1}-1)df(\nu)\nu\lrcorner\,v_{g}=0$$
Ainsi, par passage à la limite dans (\ref{eq14}), nous obtenons la formule souhaitée :
$$\lim_{r\rightarrow\infty}\int_{S_{r}}f^{-1}df(\nu)\nu\lrcorner\,v_{g}=\lim_{r\rightarrow\infty}\int_{S_{r}}df(\nu)\nu\lrcorner\,v_{g}$$
\end{proof}

\subsection{Structures de Weyl asymptotiquement plates}

Nous allons étendre la notion de platitude asymptotique aux structures de Weyl et généraliser les résultats rappelés dans la section précédente.

\begin{defn}\label{defn6}
Soient $(M,c)$ une variété conforme de dimension $n$ et $D$ une structure de Weyl sur $TM$. Nous dirons que $(M,c,D)$ est asymptotiquement plate d'ordre $\tau>0$ lorsqu'il existe une métrique $g_{0}$ dans $c$ telle que :
\begin{enumerate}
\item La variété riemannienne $(M,g_{0})$ est asymptotiquement plate d'ordre $\tau$ telle que $g_{0}\in\M_{\tau}$.
\item Pour tout entier $l$ de $1$ à $k$, la forme de Lee $\theta_{0}$ de $D$ relative à $g_{0}$
appartient à $W^{1,q}_{-\tau-1}(M^{l}_{\infty})$, où $M_{\infty}^{l}$ est le $l$-ème bout de $(M,g_{0})$ et $q>n$.
\item La codifférentielle par rapport à $g_{0}$ de la forme de Lee $\theta_{0}$ est intégrable sur $M$.
\end{enumerate}
Nous dirons que la métrique $g_{0}$ est une \emph{métrique adaptée} pour la structure de Weyl asymptotiquement plate $D$.
\end{defn}

Soient $(M,c)$ une variété conforme de dimension $n$ et $D$ une structure de Weyl sur $TM$. Supposons que $(M,c,D)$ est asymptotiquement plate à un seul bout $M_{\infty}$ et choisissons une métrique adaptée $g_{0}$ dans $c$. Sauf mention contraire, les objets riemanniens ($\nabla$, $\delta$, etc. $\ldots$) et les espaces de fonctions dans toute la suite de cette partie sont définis relativement à la métrique $g_{0}$. Nous notons $\nabla$ la connexion de Levi-Civita de $g_{0}$, $v_{0}$ sa forme volume et $\theta_{0}$ la forme de Lee de $D$ associée à $g_{0}$.

\begin{defn}\label{defn7}
Soient $(M,c,D)$ une structure de Weyl asymptotiquement plate et $g_{0}$ une métrique adaptée. Nous définissons la masse de la structure de Weyl asymptotiquement plate par :
$$m(D)(g_{0})=m(g_{0})+2(n-1)\int_{M}\delta(\theta_{0})v_{0}$$
où $m(g_{0})$ est la masse riemannienne totale de la variété $(M,g_{0})$.
La masse d'un bout $M_{\infty}^{l}$ de $M$, notée $m^{l}(D)(g_{0})$, est donnée par :
$$m^{l}(D)(g_{0})=m^{l}(g_{0})+2(n-1)\lim_{r\rightarrow\infty}\int_{S_{r}^{l}}\theta_{0}(\nu_{l})\nu_{l}\lrcorner v_{g}$$
où $m^{l}(g_{0})$ est la masse riemannienne du bout $M_{\infty}^{l}$, $S_{r}^{l}$ la sphère de rayon $r$ dans les coordonnées sur $M_{\infty}$ et $\nu_{l}$ le champ de vecteurs sortant de la sphère $S_{r}^{l}$. La masse de la variété $M$ est donc la somme des masses de chaque bout.
\end{defn}

Soit $g$ une métrique dans la classe conforme $c$ telle que $g\in M_{\tau}$; nous définissons la quantité $m(D)(g)$ par :
$$m(D)(g)=m(g)+2(n-1)\int_{M}\delta^{g}(\theta_{g})v_{g}$$
où $\delta^{g}$ et $v_{g}$ sont respectivement l'opérateur de divergence et la forme volume associés à $g$. Posons $\M_{\tau}^{0}=\{g=fg_{0}\, :\,f\in\F\}$.

\begin{prop}\label{prop20}
Soient $(M,c,D)$ une structure conforme asymptotiquement plate et $g_{0}$ une métrique adaptée.
Pour toute métrique $g$ dans $\M_{\tau}^{0}$, nous avons :
$$m(D)(g)=m(D)(g_{0})$$
\end{prop}

\begin{proof}
Soit  $g=fg_{0}$ avec $f\in\F$. Nous supposons que la variété ne possède qu'un seul bout à l'infini. Par définition, nous avons :
$$m(D)(g)=m(g)+2(n-1)\lim_{r\rightarrow\infty}\int_{S_{r}}\theta_{g}(\nu)\nu\lrcorner\,v_{g}$$
Cependant, les deux métriques sont asymptotiquement plates dans le même système de coordonnées, donc nous pouvons remplacer $v_{g}$ par $v_{0}$ dans la limite de l'intégrale. De plus, nous avons la formule :
\begin{equation}\label{eq20bis}
\theta_{g}=\theta_{0}-\frac{1}{2}\frac{df}{f}
\end{equation}
Nous obtenons :
$$m(D)(g)=m(g)+2(n-1)\lim_{r\rightarrow\infty}\int_{S_{r}}\theta_{0}(\nu)\nu\lrcorner\,v_{0}-(n-1)\lim_{r\rightarrow\infty}\int_{S_{r}}\frac{df}{f}(\nu)\nu\lrcorner\,v_{0}$$
Le théorème \ref{thm5} nous donne :
$$m(g)=m(g_{0})+(n-1)\int_{M}\Delta(f)v_{0}$$
Et, d'après le lemme \ref{lem2}, nous avons :
$$\lim_{r\rightarrow\infty}\int_{S_{r}}\frac{df}{f}(\nu)\nu\lrcorner\,v_{0}=\lim_{r\rightarrow\infty}\int_{S_{r}}df(\nu)\nu\lrcorner\,v_{0}=\int_{M}\Delta(f)v_{0}$$
Nous déduisons le résultat de ces deux dernières formules :
$$m(D)(g)=m(g_{0})+2(n-1)\int_{M}\delta(\theta_{0})v_{0}$$
\end{proof}

\subsection{Sous-classes asymptotiquement plates}

Nous savons, depuis R. Bartnik \cite{bar}, que la masse riemannienne est un invariant géométrique : la masse ne dépend pas du système de coordonnées choisi sur le bout de la variété. Plus précisément, si $g$ est asymptotiquement plate dans les coordonnées $\{z_{i}\}$ et $\{\widetilde{z}_{i}\}$, alors il existe une transformation $E$ de $\R^{n}$, composée d'une isométrie et d'une translation, telle que $\widetilde{z}=E(z)+O(r^{-\tau})$, et dans ces conditions, les masses calculées dans ces deux systèmes de coordonnées sont égales.

Soient $(M,c,D)$ une structure de Weyl asymptotiquement plate et $g_{0}$ une métrique adaptée pour $D$. Supposons que la variété ne possède qu'un seul bout à l'infini. Soit $\{z_{i}\}$ le système de coordonnées dans lequel $g_{0}$ est asymptotique à la métrique plate de $\R^{n}$. Considérons le système de coordonnées $\{\widetilde{z}_{j}\}$ sur $M_{\infty}$ tel que $\widetilde{z}_{i}=az_{i}$, où $a\in\R^{>0}$.
\begin{rem}\label{rem3}
La métrique $g_{0}$ n'est plus asymptotiquement plate dans les coordonnées $\{ \tilde{z}_{i} \}$, mais il est clair que $a^{2}g_{0}$ l'est.
\end{rem}
Notons respectivement $\widetilde{g}_{ij}$ et $g_{ij}$ les composantes de $g_{0}$ dans les coordonnées $\{\widetilde{z}_{i}\}$ et $\{z_{j}\}$, puis $\tilde{\partial}_{i}$ et $\partial_{i}$ les dérivées partielles par rapport aux coordonnées $\tilde{z}_{i}$ et $z_{i}$. Nous avons $d\tilde{z}_{i}=adz_{i}$ et $\widetilde{\partial}_{i}=a^{-1}\partial_{i}$. Par conséquent $\tilde{v}_{0}=a^{n}v_{0}$ et :
$$\widetilde{g}_{ij}=a^{-2}g_{ij}$$
En dérivant par rapport à $\widetilde{z}_{i}$, nous obtenons :
$$\widetilde{\partial}_{i}(a^{2}\widetilde{g}_{ij})=\widetilde{\partial}_{i}(g_{ij})=\sum_{k=1}^{n}\partial_{k} g_{ij}\tilde{\partial}_{i}z_{k}=a^{-1}\partial_{k} g_{ij}$$
Nous en déduisons la formule suivante :
$$\left(\widetilde{\partial}_{i}(a^{2}\widetilde{g}_{ij})-\widetilde{\partial}_{j}(a^{2}\widetilde{g}_{ii})\right)\widetilde{\partial}_{j}\lrcorner\widetilde{v_{0}}=a^{n-2}\left(\partial_{i} g_{ij}-\partial_{j} g_{ii}\right)\partial_{i}\lrcorner v_{0}$$
En intégrant sur les sphères de rayon $r$, et en faisant tendre $r$ vers l'infini, nous relions les masses de $g_{0}$ et $a^{2}g_{0}$ par :
$$m(a^{2}g_{0})=a^{n-2}m(g_{0})$$
Etudions le second membre de la masse de la structure de Weyl $D$. Rappelons que $$\int_{M}\delta(\theta_{0}) v_{0}=\int_{M}d(\theta_{0}^{\sharp}\lrcorner\,v_{0})$$
Si $g=a^{2}g_{0}$, où $a$ est une constante, nous avons $\theta_{g}=\theta_{0}$ et $\theta_{0}^{\sharp_{g}}=a^{-2}\theta_{0}^{\sharp_{0}}$, où $\sharp_{g}$ et $\sharp_{0}$ sont les isomorphismes musicaux associés à $g$ et $g_{0}$ respectivement.
Par conséquent,
$$\int_{M}\delta^{g}(\theta_{g})v_{g}=a^{n-2}\int_{M}\delta(\theta_{0})v_{0}$$
Nous venons de démontrer :

\begin{prop}\label{prop21}
Si $g_{0}$ est une métrique adaptée, nous avons :
$$m(D)(ag_{0})=a^{\frac{n-2}{2}}m(D)(g_{0})$$
\end{prop}

Nous savons désormais comment évolue notre masse lors d'un changement des coordonnées sur le bout de la variété de la forme $\widetilde{z}=aE(z)+O(r^{-\tau})$, où $E$ est une transformation euclidienne sur $\R^{n}$ et $a$ une constante positive. Nous dirons qu'un tel changement de coordonnées est asymptotiquement conforme.

Supposons désormais que notre structure de Weyl asymptotiquement plate possède $k$ bouts, où $k$ est supérieur ou égal à $1$. Dans ce cas, la masse de la variété est définie comme la somme des masses de chacun des bouts. Rappelons que $m^{l}(g_{0})$ est la masse riemannienne de $(M,g_{0})$ sur le bout $M_{\infty}^{l}$ de $M$, et $m^{l}(D)(g_{0})$ la masse de la structure de Weyl asymptotiquement plate sur le bout $M_{\infty}^{l}$de $M$. Nous avons alors :
$$m(g_{0})=\sum_{l=1}^{k}m^{l}(g_{0})$$
et
$$m^{l}(D)(g_{0})=m^{l}(g_{0})+2(n-1)\lim_{r\rightarrow\infty}\int_{S_{r}^{l}}\theta_{0}(\nu_{l})\nu_{l}\lrcorner v_{0}$$
Soit $\{a_{l}\}_{l=1\ldots k}$ une famille de $\R^{>0}$. Posons :
$$\F_{(a_{1},\ldots,a_{k})}=\{f\in C^{\infty}(M,\R^{>0})\,:\, \Delta f\in L^{1}\ et\ \forall l=1\ldots k \,,f-a_{l}\in C^{2,\alpha}_{-\tau}(M_{\infty}^{l})\}$$
En particulier, notons $\F=\F_{(1,\ldots,1)}$.

\begin{prop}\label{prop22}
Soient $g=fg_{0}$, avec $f$ dans $\F_{(a_{1},\ldots,a_{k})}$. Nous avons :
$$m(D)(g)=\sum_{l=1}^{k}a_{l}^{\frac{n-2}{2}}m^{l}(D)(g_{0})$$
En particulier, si $g_{1}=f_{1}g_{0}$ et $g_{2}=f_{2}g_{0}$, avec $f_{1}$ et $f_{2}$ dans $\F_{(a_{1},\ldots,a_{k})}$, nous avons :
$$m(D)(g_{1})=m(D)(g_{2})$$
\end{prop}

\begin{proof}
Fixons un bout $M_{\infty}^{l}$ de $(M,g_{0})$. Sur $M_{\infty}^{l}$, nous avons $f=a_{l}f_{l}$, où $f_{l}$ appartient à $\F$ sur $M_{\infty}^{l}$. D'après le théorème \ref{thm5}, la métrique $g_{l}=f_{l}g_{0}$ est asymptotiquement plate sur $M_{\infty}^{l}$ et la proposition \ref{prop20}, adaptée pour le bout $M_{\infty}^{l}$, donne :
$$m^{l}(D)(g_{l})=m^{l}(D)(g_{0})$$
Pour la métrique adaptée $g_{l}$, en utilisant la proposition \ref{prop21}, nous obtenons :
$$m^{l}(D)(g)=a_{l}^{\frac{n-2}{2}}m^{l}(D)(g_{l})$$
Donc, nous avons :
$$m(D)(g)=\sum_{l=1}^{k}a_{l}^{\frac{n-2}{2}}m^{l}(D)(g_{0})$$
\end{proof}

Par conséquent, la masse de $D$ évaluée en $g=fg_{0}$ ne dépend pas de la fonction $f$ choisie dans $\F_{(a_{1},\ldots,a_{k})}$, mais uniquement des réels $a_{l}$ et de la métrique de référence $g_{0}$.
Une métrique adaptée étant fixée, nous remarquons que les ensembles de fonctions $\F_{a_{1},\ldots,a_{k}}$ définissent des sous-classes de métriques de la classe conforme $c$ pour lesquelles la masse de la structure de Weyl est invariante. Notons $\Gamma=\{\F_{a_{1},\ldots,a_{k}}:(a_{1},\ldots,a_{n})\in(\R^{>0})^{n}\}$. Nous avons alors une application :
$$m(D)\ :\  \Gamma\rightarrow\R$$
définie par :
$$m(D)(\F_{a_{1},\ldots,a_{k}})=m(D)(fg_{0})$$
où $f$ est une fonction quelconque de $\F_{a_{1},\ldots,a_{k}}$. L'espoir de cette remarque est de parvenir à démontrer que $\Gamma$ constitue l'ensemble de toutes les métriques asymptotiquement plates de la classe conforme $c$. Dans ce cas, nous aurons défini une masse conforme $m(D)$ indépendante du choix de la métrique adaptée.

\subsection{Théorème de la masse conforme positive}

Soient $(M,c,D)$ une structure de Weyl asymptotiquement plate ne possédant qu'un seul bout $M_{\infty}$ et $g_{0}$ une métrique adaptée pour $D$. Soit $\nabla$ la connexion de Levi-Civita de $g_{0}$. Supposons que $M$ est une variété spinorielle. Notons $\Sigma=\Sigma^{(\frac{2-n}{2})}$ l'espace des spineurs conformes de poids $(2-n)/2$. Nous notons désormais $\Sigma^{0}$ le fibré des spineurs riemanniens relatif à la métrique $g_{0}$ et $\D^{0}$ l'opérateur de Dirac induit par $\nabla$ sur $\Sigma^{0}$. Rappelons que les espaces de fonctions considérés sont relatifs à la métrique adaptée $g_{0}$.

\begin{lem}\label{lem3}
Soit $\psi$ une section du fibré des spineurs $\Sigma^{0}$. Supposons que $\psi$ est $D$-parallèle, $i.e$ $D\psi=0$. Si $\psi$ tend vers $0$ dans $M_{\infty}$, $i.e.$ $\lim_{|x|\rightarrow\infty}|\psi(x)|=0$, alors $\psi$ est identiquement nulle.
\end{lem}

\begin{proof}
Soient $x$ dans $M$ et $\{e_{i}\}$ une base de $T_{x}M$. Calculons la différentielle de $|\psi|^{2}$ :  $d|\psi|^{2}=\sum_{i=1}^{n}\nabla_{e_{i}}(|\psi|^{2})e_{i}^{\ast}$. La compatibilité de la connexion de Levi-Civita avec le produit scalaire donne $\nabla_{e_{i}}(|\psi|^{2})=(\nabla_{e_{i}}\psi,\psi)+(\psi,\nabla_{e_{i}}\psi)$. Comme $\psi$ est $D$-parallèle, la formule (\ref{eq3ter}) donne : $$\nabla_{e_{i}}\psi=\frac{n-1}{2}\theta_{0}(e_{i})\psi+\frac{1}{2}e_{i}\cdot\theta_{0}\cdot\psi$$
Nous obtenons :
\begin{eqnarray}
\nabla_{e_{i}}(|\psi|^{2})&=&(n-1)\theta_{0}(e_{i})|\psi|^{2}+\frac{1}{2}((e_{i}\cdot\theta_{0}+\theta_{0}\cdot e_{i})\cdot\psi,\psi)\nonumber\\
&=&(n-2)\theta_{0}(e_{i})|\psi|^{2}\nonumber
\end{eqnarray}
Par conséquent, $d|\psi|^{2}=(n-2)|\psi|^{2}\theta_{0}$. Soient $x_{0}$ dans $M_{\infty}$ et $\gamma$ une géodésique paramétrée sur $[0,+\infty[$ telle que $\gamma(0)=x_{0}$ et  $\lim_{t\rightarrow\infty}|\gamma(t)|=+\infty$. Posons $\phi(t)=|\psi|^{2}_{\gamma(t)}$ et $f(t)=\theta_{0}(\dot{\gamma}(t))$. La fonction $\phi$ est donc solution de l'équation différentielle ordinaire $\phi^{\prime}(t)=(n-2)f(t)\phi(t)$ sur $[0,\infty[$. Par conséquent, pour tout $t$, $\phi(t)=c\exp(F(t))$, où $c$ est une constante et $F$ une primitive de $f$. Cependant, par hypothèse $\lim_{t\rightarrow\infty}\phi(t)=0$, donc $\phi$ est identiquement nulle. Ceci étant vrai pour toute courbe allant vers l'infini sur $M_{\infty}$, nous en déduisons que $\psi$ est identiquement nul.
\end{proof}

Soit $\psi$ une section de $\Sigma^{0}$. Soit $\omega$ une $1$-forme sur $M$. La formule de Stokes nous donne
$$-\int_{M_{r}}\delta(\omega)v_{0}=\int_{S_{r}}\omega(\nu)\nu\lrcorner\,v_{0}$$
où $\nu$ est le champ de vecteurs normal unitaire sortant de $S_{r}$, et rappelons que $M_{r}$ est le compact de $M$ défini par $M_{r}=M\setminus E_{r}$ dont le bord s'identifie à la sphère $S_{r}$ de $\R^{n}$. Par conséquent, si la limite existe, nous avons : 
$$-\int_{M}\delta(\omega)v_{0}=\lim_{r\rightarrow\infty}\int_{S_{r}}\omega(\nu)\nu\lrcorner\,v_{g}$$

\begin{prop}\label{prop23}
Soit $1-n<\delta<0$. Supposons que la courbure scalaire $Scal^{D}$ de la connexion de Weyl est positive. Dans ces conditions, l'opérateur de Dirac de poids conforme $(2-n)/2$, $\D^{2}$ : $W^{2,q}_{-\tau}(\Sigma)\rightarrow L^{q}_{-\tau-2}(\Sigma)$, est un isomorphisme.
\end{prop}

\begin{proof}
Montrons que l'opérateur $\D^{2}$ est bien défini sur ces espaces de Sobolev. Soit $\psi$ dans $W^{2,q}_{-\tau}(\Sigma)$.
La proposition \ref{prop4} et l'inégalité triangulaire nous donnent : 
$$|\D\psi|\leq|\D^{0}\psi|+\frac{1}{2}|\theta_{0}||\psi|$$
Par hypothèse sur $\psi$, la norme $|\psi|$ est bornée sur $M$ (puisque $|\psi|=O(r^{-\tau})$) et la proposition \ref{prop18} nous donne $\D^{0}\psi\in W_{-\tau-1}^{1,q}$. De plus, la condition ($2$) de la définition \ref{defn6} assure que $\theta_{0}$ appartient également à $W_{-\tau-1}^{1,q}$. L'inégalité ci-dessus permet de conclure que $\D\psi$ est une section de $W_{-\tau-1}^{1,q}$. Par le même raisonnement en remplaçant  $\psi$ par $\D\psi\in W^{1,q}_{-\tau-1}$, nous démontrons que $\D^{2}\psi$ appartient à $W^{0,q}_{-\tau-2}$. Soit $\psi$ dans $W^{2,q}_{-\tau}$ tel que $\D^{2}\psi=0$, montrons que $\psi$ est identiquement nul. Le corollaire \ref{cor13} nous donne :
$$|D\psi|_{h}^{2}+\frac{1}{4}Scal^{D}|\psi|_{h}^{2}=-\delta\big(h(\psi,D\psi)\big)$$
En intégrant cette formule sur le compact $M_{r}$ de $M$, pour $r>R$, nous avons : 
$$\int_{M_{r}}|D\psi|_{h}^{2}+\frac{1}{4}\int_{M_{r}}Scal^{D}|\psi|_{h}^{2}=-\int_{M_{r}}\delta\big(h(\psi,D\psi)\big)$$
La métrique de référence $g_{0}$ étant fixé, la divergence conforme s'identifie avec la divergence riemannienne relative à $g_{0}$. Donc la formule de Stokes nous donne :
\begin{equation}\label{eq15}
\int_{M_{r}}|D\psi|_{h}^{2}+\frac{1}{4}\int_{M_{r}}Scal^{D}|\psi|_{h}^{2}=\int_{S_{r}}(\psi,D_{\nu}\psi)\nu\lrcorner v_{0}
\end{equation}
où $\nu$ est le champ de vecteurs unitaire sortant de $S_{r}$.
La formule (\ref{eq3ter}) reliant la structure de Weyl et la connexion de Levi-Civita de la métrique $g_{0}$ nous donne :  
$$(\psi,D_{\nu}\psi)=(\psi,\nabla_{\nu}\psi)+\frac{1-n}{2}\theta_{0}(\nu)|\psi|^{2}+\frac{1}{2}(\nu\cdot\psi,\theta_{0}\cdot\psi)$$
Par l'inégalité triangulaire et de Cauchy-Schwarz, nous obtenons :
$$|(\psi,D_{\nu}\psi)|\leq|\nabla\psi||\psi|+\frac{n-1}{2}|\theta_{0}||\psi|^{2}+\frac{1}{2}|\theta_{0}||\psi|^{2}$$
Par hypothèse sur la section $\psi$ et la forme de Lee $\theta_{0}$, nous avons  $(\psi,D_{\nu}\psi)=O(r^{-2\tau-1})$, avec $2\tau+1>n-1$. Par conséquent, nous arrivons à l'égalité suivante :
$$\lim_{r\rightarrow\infty}\int_{S_{r}}(\psi,D_{\nu}\psi)\nu\lrcorner\,v_{0}=0$$
Nous en déduisons, par passage à la limite quand $r$ tend vers l'infini dans (\ref{eq15}), la formule suivante : 
$$\int_{M}|D\psi|_{h}^{2}+\frac{1}{4}\int_{M}Scal^{D}|\psi|_{h}^{2}=0$$
La courbure scalaire $Scal^{D}$ de la connexion de Weyl est positive, donc les deux termes intégrés sont nuls. En particulier, $D\psi=0$ et $\psi$ tend vers $0$ à l'infini, donc, d'après le lemme \ref{lem3}, $\psi$ est identiquement nul. Nous avons donc démontré que l'opérateur $\D^{2}$ : $W_{-\tau}^{2,q}\rightarrow L_{-\tau-2}^{q}$ est injectif. Cependant, la proposition \ref{prop14} démontre que $\D^{2}$ est formellement autoadjoint : 
$$\D^{2}=(\D^{2})^{\ast}\,: \,W_{\tau+2-n}^{2,q^{\prime}}\rightarrow L^{q^{\prime}}_{\tau-n}$$
Un raisonnement similaire démontre que cet opérateur est injectif (par hypothèse nous avons $1-n<\tau+2-n<0$).  L'opérateur $\D^{2}$ est de Fredholm de noyau et conoyau trivial donc $\D^{2}$ est un isomorphisme de $W^{2,q}_{-\tau}$ dans $L^{q}_{-\tau-2}$ (voir \cite{bar}).
\end{proof}

Sous les hypothèses de cette section, nous sommes en mesure de démontrer la convergence vers la masse conforme de l'intégrale du terme de divergence dans la formule de Lichnerowicz conforme II.

\begin{prop}\label{prop24}
Soit $\psi$ un spineur asymptotiquement constant. La métrique adaptée $g_{0}$ étant fixée, nous avons la formule suivante :
$$-\int_{M}\delta(\omega_{\psi})=\frac{1}{4}m(D)(g_{0})|\psi_{0}|^{2}$$
où $\psi_{0}$ est le spineur constant tel que $\psi-\psi_{0}\in W_{-\tau}^{2,q}$, et où la section $\omega_{\psi}$ de $T^{\ast}M\otimes L^{(2-n)}$ est définie par : $\omega_{\psi}(X)=h\big(\psi,X^{\flat}\cdot\D\psi+D_{X}\psi\big)$.
\end{prop}

\begin{proof}
La métrique $g_{0}$ étant donnée, la proposition \ref{prop4} et la formule (\ref{eq3ter}) appliquées pour $g_{0}$ permettent d'identifier $\omega_{\psi}$ à une $1$-forme sur $M$ par la formule suivante :
\begin{equation}\label{eq16}
\omega_{\psi}(X)=\omega^{0}_{\psi}(X)+\frac{1-n}{2}\theta_{0}(X)|\psi|^{2}
\end{equation}
avec $\omega_{\psi}^{0}(x)=(\psi,X\cdot\D^{0}\psi+\nabla_{X}\psi)$. Soit $\psi_{0}$ le spineur constant tel que $\psi-\psi_{0}\in W_{-\tau}^{2,q}$; la proposition \ref{prop19} donne :
\begin{equation}\label{eq24}
-\int_{M}\delta(\omega_{\psi}^{0})v_{0}=\lim_{r\rightarrow\infty}\int_{S_{r}}\omega_{\psi}^{0}(\nu)\nu\lrcorner v_{0}=\frac{1}{4}m(g_{0})|\psi_{0}|^{2}
\end{equation}
où $\nu$ est le champ de vecteurs unitaire sortant de $S_{r}$. Posons $\psi_{1}=\psi-\psi_{0}$, nous avons : $$|\psi|^{2}=|\psi_{1}|^{2}+|\psi_{0}|^{2}+2\Re\big((\psi_{1},\psi_{0})\big)$$
Cependant, $\psi_{1}\in W_{-\tau}^{2,q}(\Sigma)$ et $\theta_{0}\in W_{-\tau-1}^{1,q}$. Donc :  
$$\theta_{0}(\nu)\left(|\psi_{1}|^{2}+2\Re\big((\psi_{1},\psi_{0})\big)\right)=O(r^{-2\tau-1})$$
avec $2\tau+1>n-1$. Par conséquent, nous obtenons l'égalité suivante :
\begin{equation}\label{eq24bis}
\lim_{r\rightarrow\infty}\int_{S_{r}}|\psi|^{2}\theta_{0}(\nu)\nu\lrcorner v_{0}=|\psi_{0}|^{2}\lim_{r\rightarrow\infty}\int_{S_{r}}\theta_{0}(\nu)\nu\lrcorner v_{0}
\end{equation}
Remarquons que cette limite existe puisque la codifférentielle de la forme de Lee $\theta_{0}$ est intégrable (condition ($3$) de la définition \ref{defn6}) et donc (\ref{eq24bis}) devient :
\begin{equation}\label{eq24ter}
\int_{M}\delta(|\psi|^{2}\theta_{0})v_{0}=\big(\int_{M}\delta(\theta_{0})v_{0}\big)|\psi_{0}|^{2}
\end{equation}
Les équations (\ref{eq24}) et (\ref{eq24ter}) nous donnent l'égalité souhaitée :
$$-\int_{M}\delta(\omega_{\psi})v_{0}=\frac{1}{4}\big(m(g_{0})+2(n-1)\int_{M}\delta(\theta_{0})v_{0}\big)|\psi_{0}|^{2}$$
\end{proof}

\begin{thm}\label{thm6}
Soient $(M,c)$ une variété conforme spinorielle et $D$ une structure de Weyl asymptotiquement plate sur $M$. Supposons que la courbure scalaire $Scal^{D}$ de $D$ est positive. Alors la masse $m(D)$ associée à $D$ est positive et cette masse est nulle si et seulement si $(M,c)$ est isomorphe à l'espace $\R^{n}$ munit de la classe conforme canonique.
\end{thm}

\begin{proof}
Soient $g_{0}$ une métrique adaptée pour $D$ et $\psi_{0}$ un spineur constant. Les hypothèses de décroissance asymptotique sur $D$ nous donnent $\D\psi_{0}\in W^{1,q}_{-\tau-1}\cap L^{q}_{-\tau-2}$. D'après la proposition \ref{prop23}, il existe $\psi_{1}\in W_{-\tau}^{2,q}$ tel que $\D^{2}\psi_{1}=-\D\psi_{0}$. Cependant, par régularité elliptique, $\psi_{1}\in W^{3,q}_{-\tau+1}$ et donc $\psi_{2}:=\D\psi_{1}\in W^{2,q}_{-\tau}$. Posons $\psi=\psi_{2}+\psi_{0}$. La section $\psi$ est $\D$-harmonique et asymptotique au spineur constant $\psi_{0}$. En intégrant sur $M_{r}$ la formule de Lichnerowicz conforme II (théorème \ref{thm3}), nous obtenons :
$$\int_{M_{r}}|D\psi|^{2}v_{0}+\frac{1}{4}\int_{M_{r}}Scal^{D}|\psi|^{2}v_{0}=-\int_{M_{r}}\delta(\omega_{\psi})v_{0}$$
La proposition \ref{prop24}, par passage à la limite, nous donne alors :
$$\int_{M}|D\psi|^{2}v_{0}+\frac{1}{4}\int_{M}Scal^{D}|\psi|^{2}v_{0}=m(D)(g_{0})|\psi_{0}|^{2}$$
Cependant, la courbure scalaire $Scal^{D}$ de la structure de Weyl $D$ est positive, donc la masse $m(D)$ associée à $D$ est positive d'après la formule précédente.

Supposons que la masse $m(D)$ soit nulle. Nous avons alors un spineur $\psi$ dans $\Sigma$ tel que :
$$\int_{M}|D\psi|_{h}^{2}+\frac{1}{4}\int_{M}Scal^{D}|\psi|_{h}^{2}=0$$
Les deux termes de cette équation sont positifs, donc nuls. La section $\psi$ est alors $D$-parallèle. Par conséquent, $|\psi|_{h}^{\frac{2}{2-n}}$ est une section $\nabla^{D}$-parallèle du fibré $L$. En effet, le calcul donne :
$$\nabla^{D}|\psi|_{h}^{\frac{2}{2-n}}=\nabla^{D}\big(|\psi|_{h}^{2}\big)^{\frac{1}{2-n}}=\frac{1}{2-n}|\psi|_{h}^{\frac{n-1}{2-n}}\big(h(D\psi,\psi)+h(\psi,D\psi)\big)$$
Ainsi, la structure de Weyl $D$ est exacte, $i.e.$ il existe une métrique $g$ dans la classe conforme $c$ dont $D$ est la connexion de Levi-Civita. \'Ecrivons $g=fg_{0}$, où $f$ est une fonction strictement positive sur $M$. La structure de Weyl $D$ étant la connexion de Levi-Civita de $g$, la forme de Lee $\theta_{g}$ associée est identiquement nulle et la formule (\ref{eq20bis}) nous donne :
$$\theta_{0}=\frac{1}{2}\frac{df}{f}$$
Posons $h=\log(f)$. Nous avons $dh=\frac{df}{f}$, et $h\in C^{1,\alpha}_{-\tau-1}(M)$ par hypothèse sur $\theta_{0}$.
\begin{lem}\label{lem4}
Soit $h$ une fonction $C^{\infty}$ sur $M$. Supposons que $dh$ appartient à $C^{1,\alpha}_{-\tau-1}(M)$. La fonction $h$ possède une limite finie $a$ telle que $h-a\in C^{2,\alpha}_{-\tau}(M)$.
\end{lem}
\begin{proof}
Soient $x$ et $y$ dans $E_{R}=\R^{n}\setminus B_{R}$ tels que $|y|\geq|x|\geq R_{1}$, pour $R_{1}>R$. Soit $z=\frac{|x|}{|y|}y$ dans $\R^{n}$. Considérons $X_{t}$ le grand arc de cercle paramétré sur $[0,1]$ joignant $x$ et $z$. Nous avons $|X_{t}|=|x|\geq R_{1}$ et $|\dot{X}_{t}|\leq\pi|x|$. L'égalité des accroissement finis entre $0$ et $1$ (pour $f(t)=h(X_{t})$) nous donne :
$$h(x)-h(z)=\int_{0}^{1}dh_{X_{t}}(\dot{X}_{t})dt$$
Par hypothèse, $|dh_{X_{t}}|=O(|X_{t}|^{-\tau-1})$, $|X_{t}|=|x|$ et $|\dot{X}_{t}|\leq\pi|x|$, nous obtenons alors l'inégalité suivante :
$$|h(x)-h(z)|\leq \pi C|x|^{-\tau}$$
où $C$ est une constante positive indépendante de $x$ et $y$. Comme $|x|\geq R_{1}$, cette inégalité implique :
\begin{equation}\label{eqlem4a}
|h(x)-h(z)|\leq \pi CR_{1}^{-\tau}
\end{equation}
Posons $T=\frac{|x|}{|y|}$; l'égalité des accroissement finis entre $z$ et $y$ nous donne :
$$h(Ty)-h(y)=\int_{0}^{1}dh_{y+t(Ty-y)}(Ty-y)dt$$
Par le même raisonnement que précédemment, nous obtenons l'inégalité suivante :
$$|h(Ty)-h(y)|\leq C|y|^{-\tau}\int_{0}^{1}(1-T)\big(1+t(T-1)\big)^{-\tau-1}dt$$
Nous pouvons calculer l'intégrale :
$$\int_{0}^{1}(1-T)\big(1+t(T-1)\big)^{-\tau-1}dt=\tau^{-1}\left(\frac{|x|^{-\tau}}{|y|^{-\tau}}-1\right)$$
Nous en déduisons :
\begin{equation}\label{eqlem4b}
|h(Ty)-h(y)|\leq\tau^{-1}C(|x|^{-\tau}-|y|^{-\tau})\leq\tau^{-1}CR_{1}^{-\tau}
\end{equation}
Par conséquent, (\ref{eqlem4a}) et (\ref{eqlem4b}) établissent l'existence d'une constante $C_{1}$ indépendante de $x$ et $y$ telle que :
$$|h(x)-h(y)|\leq C_{1}R_{1}^{-\tau}$$
Ainsi, pour tout $\varepsilon>0$, nous pouvons choisir $R_{1}$ suffisamment grand tel que, pour tout $x$ et $y$ dans $E_{r}$ vérifiant $|y|\geq|x|\geq R_{1}$, $|h(x)-h(y)|\leq\varepsilon$. Par conséquent, d'après le critère de Cauchy, $h(x)$ admet une limite finie lorsque $x$ tend vers l'infini. Par passage à la limite quand $T$ tend vers l'infini dans l'inégalité (\ref{eqlem4b}), nous obtenons :
\begin{equation}\label{eqlem4}
|h(x)-a|\leq \frac{C}{\tau}|x|^{-\tau}
\end{equation}
pour tout $x$ dans $E_{R}$. En observant la norme des espaces de Hölder à poids, nous pouvons remarquer que l'estimation (\ref{eqlem4}) et l'hypothèse $dh\in C_{-\tau-1}^{1,\alpha}$ suffisent pour démontrer que $h$ appartient à $C_{-\tau}^{2,\alpha}$.
\end{proof}
En appliquant le lemme \ref{lem4} à la fonction $h=\log(f)$, nous démontrons que $f$ possède une limite finie strictement positive $b=\exp(a)$ à l'infini et que $f-b\in C^{2,\alpha}_{-\tau}$. En changeant $f$ par $b^{-1}f$, nous ne modifions pas la connexion de Levi-Civita associée à $g$. Nous pouvons donc supposer que $g=fg_{0}$ avec $f\in\F$. D'après le théorème \ref{thm5}, la métrique $g$ est asymptotiquement plate. Nous avons donc $\psi\in\Sigma^{g}$ tel que $D\psi=0$, où $D$ est la connexion de Levi-Civita de la métrique asymptotiquement plate $g$, la preuve du théorème de la masse positive \cite{bar} démontre dans ce cas que $M$ est isométrique à $\R^{n}$ munit de la métrique plate. En conclusion, la variété conforme $(M,c)$ est isomorphe à l'espace $\R^{n}$ munit de sa classe conforme canonique.
\end{proof}


\labelsep .5cm

 \end{document}